\documentclass{article}
\usepackage{amssymb}

\input{tcilatex}
\begin{document}

\begin{center}
\textbf{Unique continuation properties for Schredinger operators in Hilbert
spaces}

\textbf{Veli Shakhmurov}

Department of Mechanical Engineering, Okan University, Akfirat, Tuzla 34959
Istanbul, Turkey,

E-mail: veli.sahmurov@okan.edu.tr

Institute of Mathematics and Mechanics, Azerbaijan National Academy of
Sciences, Azerbaijan, AZ1141, Baku, F. Agaev, 9,

E-mail: veli.sahmurov@gmail.com

\bigskip

A\textbf{bstract}
\end{center}

Here, the Morgan type uncertainty principle and unique continuation
properties of abstract Schredinger equations with time dependent potentials
are obtained in Hilbert space valued function\ classes. The equations
include linear operator in abstract Hilbert spaces $H$ dependent on space
variables. So, by selecting appropriate spaces $H$ and operators, we derive
unique continuation properties for numerous classes of Schr\"{o}dinger type
equations and its systems, which occur in a wide variety of physical systems.

\textbf{Key Word:}$\mathbb{\ \ }$Schredinger equations\textbf{, }Positive
operators\textbf{, }Semigroups of operators, Unique continuation, Morgan
type uncertainty principle

\begin{center}
\bigskip\ \ \textbf{AMS 2010: 35Q41, 35K15, 47B25, 47Dxx, 46E40 }

\textbf{1. Introduction, definitions}
\end{center}

\bigskip In this paper, the unique continuation properties of the abstract
Schr\"{o}dinger equations

\begin{equation}
i\partial _{t}u+\Delta u+A\left( x\right) u+V\left( x,t\right) u=0,\text{ }%
x\in R^{n},\text{ }t\in \left[ 0,1\right] ,  \tag{1.1}
\end{equation}%
are studied, where $A=A\left( x\right) $ is a linear operator$,$ $V\left(
x,t\right) $ is a given potential operator function in a Hilbert space $H$,
the subscript $t$ indicates the partial derivative with respect to $t$, $%
\Delta $ denotes the Laplace operator in $R^{n}$ and $u=$ $u(x,t)$ is the $H$%
-valued unknown function. The goal is to obtain sufficient conditions on $A$%
, the potential $V$ and the behavior of the solution $u$ at two different
times, $t_{0}=0$ and $t_{1}=1,$ which guarantee that $u\equiv 0$ in $%
R^{n}\times \lbrack 0,1]$. This linear result is then applied to show that
two regular solutions $u_{1}$ and $u_{2}$ of non-linear Schredinger equation 
\begin{equation}
i\partial _{t}u+\Delta u+Au=F\left( u,\bar{u}\right) ,\text{ }x\in R^{n},%
\text{ }t\in \left[ 0,1\right] ,  \tag{1.2}
\end{equation}%
for general non-linearities $F$, must agree in $R^{n}\times \lbrack 0,1]$,
when $u_{1}-u_{2}$ and its gradient decay faster than any quadratic
exponential at times $0$ and $1$.

Unique continuation properties for Schredinger equations have been studied
by $\left[ \text{5-8}\right] ,$ $\left[ \text{21-23}\right] $ and the
references therein.\ In contrast to the mentioned above results we will
study the unique continuation properties of abstract Schredinger equations
with operator potentials. Abstract differential equations were studied e.g.
in $\left[ 1\right] $, $\left[ 10\right] ,$ $\left[ 13-19\right] $, $\left[
25\text{, }26\right] .$ Since the Hilbert space $H$ is arbitrary and $A$, $V$
are possible linear operators, by choosing $H$, $A$ and $V$\ we can obtain
numerous classes of Schr\"{o}dinger type equations and its systems which
occur in different applications. If we\ choose the abstract space $H$ a
concrete Hilbert space, for example $H=L^{2}\left( G\right) $, $A=L,$ $%
V=q\left( x,t\right) $, where $G$ is a domain in $R^{m}$ with smooth
boundary, $L$ is a elliptic operator with respect to the variable $y\in G$
and $q$ is a complex valued function, then we obtain the unique continuation
properties of the following Schr\"{o}dinger equation%
\begin{equation}
\partial _{t}u=i\left[ \Delta u+Lu+q\left( x,t\right) u\right] ,\text{ }x\in
R^{n},\text{ }y\in G,\text{ }t\in \left[ 0,1\right] ,  \tag{1.3}
\end{equation}%
where $u=u\left( x,y,t\right) $. \ Moreover, let $H=L^{2}\left( 0,1\right) $
and let $A$ to be a differential operator with generalized Wentzell-Robin
boundary condition defined by 
\begin{equation}
D\left( A\right) =\left\{ u\in W^{2,2}\left( 0,1\right) ,\text{ }%
B_{j}u=Au\left( j\right) =0\text{, }j=0,1\right\} ,\text{ }  \tag{1.4}
\end{equation}%
\[
\text{ }Au=au^{\left( 2\right) }+bu^{\left( 1\right) }, 
\]%
where $a$ and $b$ are real valued functions on $\ R^{n}\times \left[ 0,1%
\right] $. Then, we obtain the unique continuation properties of the
Wentzell-Robin type boundary value problem (BVP) for the nonlinear
Schredinger type equation 
\begin{equation}
i\partial _{t}u+\Delta u+a\frac{\partial ^{2}u}{\partial y^{2}}+b\frac{%
\partial u}{\partial y}=F\left( u,\bar{u}\right) ,\text{ }  \tag{1.5}
\end{equation}%
\[
u=u\left( x,y,t\right) \text{, }x\in R^{n},\text{ }y\in G,\text{ }t\in \left[
0,1\right] , 
\]%
\ \ \ 

\begin{equation}
a\left( j\right) u_{yy}\left( x,j,t\right) +b\left( j\right) u_{y}\left(
x,j,t\right) =0\text{, }j=0,1\text{, for a.e. }x\in R^{n},\text{ }t\in
\left( 0,1\right) .  \tag{1.6}
\end{equation}

Note that, the regularity properties of Wentzell-Robin type BVP for elliptic
equations have been studied e.g. in $\left[ \text{11, 12 }\right] $ and the
references therein. Moreover, if put $H=l_{2}$ and choose $A$ and $V\left(
x,t\right) $ as infinite matrices $\left[ a_{mj}\right] $, $\left[
b_{mj}\left( x,t\right) \right] ,$ $m,$ $j=1,2,...,\infty $ respectively,
then we obtain the unique continuation properties of the following system of
Schredinger equation 
\begin{equation}
\partial _{t}u_{m}=i\left[ \Delta u_{m}+\sum\limits_{j=1}^{N}\left(
a_{mj}+b_{mj}\left( x,t\right) \right) u_{j}\right] ,\text{ }x\in R^{n},%
\text{ }t\in \left( 0,1\right) .  \tag{1.7}
\end{equation}

Let $E$ be a Banach space and $\gamma =\gamma \left( x\right) $ be a
positive measurable function on a domain $\Omega \subset R^{n}.$ Here, $%
L_{p,\gamma }\left( \Omega ;E\right) $ denotes the space of strongly
measurable $E$-valued functions that are defined on $\Omega $ with the norm

\[
\left\Vert f\right\Vert _{_{p,\gamma }}=\left\Vert f\right\Vert
_{L_{p,\gamma }\left( \Omega ;E\right) }=\left( \int \left\Vert f\left(
x\right) \right\Vert _{E}^{p}\gamma \left( x\right) dx\right) ^{\frac{1}{p}},%
\text{ }1\leq p<\infty , 
\]

\[
\left\Vert f\right\Vert _{L_{\infty ,\gamma }\left( \Omega ;E\right)
}=ess\sup\limits_{x\in \Omega }\left\Vert f\left( x\right) \right\Vert
_{E}\gamma \left( x\right) ,\text{ }p=\infty . 
\]

For $\gamma \left( x\right) \equiv 1$ the space $L_{p,\gamma }\left( \Omega
;E\right) $ will be denoted by $L^{p}=L^{p}\left( \Omega ;E\right) $ for $%
p\in \left[ 1,\infty \right] .$

Let $\mathbf{p=}\left( p_{1},p_{2}\right) $ and $\Omega =\Omega _{1}\times
\Omega _{2}$, where $\Omega _{k}\in R^{n_{k}}$. $L_{x}^{p_{1}}L_{t}^{p_{2}}%
\left( \Omega ;E\right) $ will denote the space of all $E$-valed, $\mathbf{p}
$-summable\ functions with mixed norm, i.e., the space of all measurable
functions $f$ defined on $G$ equipped with norm%
\[
\left\Vert f\right\Vert _{L_{x}^{p_{1}}L_{y}^{p_{2}}\left( \Omega ;E\right)
}=\left( \left( \dint\limits_{\Omega _{1}}\left\Vert f\left( x,t\right)
\right\Vert _{E}^{p_{2}}dt\right) ^{\frac{p_{1}}{p_{2}}}dx\right) ^{\frac{1}{%
p_{1}}}. 
\]

Let $H$ be a Hilbert space and 
\[
\left\Vert u\right\Vert =\left\Vert u\right\Vert _{H}=\left( u,u\right)
_{H}^{\frac{1}{2}}=\left( u,u\right) ^{\frac{1}{2}}\text{ for }u\in H. 
\]%
\ For $p=2$ and $E=H$ is a Hilbert space, $L^{p}\left( \Omega ;E\right) $ is
to be a Hilbert space with inner product 
\[
\left( f,g\right) _{L^{2}\left( \Omega ;H\right) }=\int\limits_{\Omega
}\left( f\left( x\right) ,g\left( x\right) \right) _{H}dx 
\]%
for $f$, $g\in L^{2}\left( \Omega ;H\right) $.

Let $C\left( \Omega ;E\right) $ denote the space of $E-$valued, bounded
uniformly continuous functions on $\Omega $ with norm 
\[
\left\Vert u\right\Vert _{C\left( \Omega ;E\right) }=\sup\limits_{x\in
\Omega }\left\Vert u\left( x\right) \right\Vert _{E}. 
\]

$C^{m}\left( \Omega ;E\right) $\ will denote the space of $E$-valued bounded
uniformly strongly continuous and $m$-times continuously differentiable
functions on $\Omega $ with norm 
\[
\left\Vert u\right\Vert _{C^{m}\left( \Omega ;E\right) }=\max\limits_{0\leq
\left\vert \alpha \right\vert \leq m}\sup\limits_{x\in \Omega }\left\Vert
D^{\alpha }u\left( x\right) \right\Vert _{E}. 
\]

Moreover, $C_{0}^{\infty }\left( \Omega ;E\right) -$denotes the space of $E$%
-valued infinity many differentiable finite functions.

Let

\[
O_{r}=\left\{ x\in R^{n},\text{ }\left\vert x\right\vert <r\right\} \text{, }%
r>0. 
\]%
Let $\mathbb{N}$ denote the set of all natural numbers, $\mathbb{C}$ denote
the set of all complex numbers.

Let $\Omega $ be a domain in $R^{n}$ and let $m$ be a positive integer$.$\ $%
W^{m,p}\left( \Omega ;E\right) $ denotes the space of all functions $u\in
L^{p}\left( \Omega ;E\right) $ that have generalized derivatives $\frac{%
\partial ^{m}u}{\partial x_{k}^{m}}\in L^{p}\left( \Omega ;E\right) ,$ $%
1\leq p\leq \infty $ with the norm 
\[
\ \left\Vert u\right\Vert _{W^{m,p}\left( \Omega ;E\right) }=\left\Vert
u\right\Vert _{L^{p}\left( \Omega ;E\right)
}+\sum\limits_{k=1}^{n}\left\Vert \frac{\partial ^{m}u}{\partial x_{k}^{m}}%
\right\Vert _{L^{p}\left( \Omega ;E\right) }<\infty . 
\]

Let $E_{0}$ and $E$ be two Banach spaces and suppose $E_{0}$ is continuously
and densely embedded into $E$. Here,\ $W^{m,p}\left( \Omega ;E_{0},E\right) $
denotes the space $W^{m,p}\left( \Omega ;E\right) \cap $ $L^{p}\left( \Omega
;E\right) $ equipped with norm 
\[
\ \left\Vert u\right\Vert _{W^{m,p}\left( \Omega ;E_{0,}E\right)
}=\left\Vert u\right\Vert _{L^{p}\left( \Omega ;E_{0}\right)
}+\sum\limits_{k=1}^{n}\left\Vert \frac{\partial ^{m}u}{\partial x_{k}^{m}}%
\right\Vert _{L^{p}\left( \Omega ;E\right) }<\infty . 
\]

Let $E_{1}$ and $E_{2}$ be two Banach spaces. Let $L\left(
E_{1},E_{2}\right) $ will denote the space of all bounded linear operators
from $E_{1}$ to $E_{2}.$ For $E_{1}=E_{2}=E$ it will be denoted by $L\left(
E\right) .$

Let $A$ is a symmetric operator in a Hilbert sapace $H$ with domain $D\left(
A\right) .$ Here, $H\left( A\right) $ denotes the domain of $A$ equipped
with graphical norm, i.e. 
\[
\left\Vert u\right\Vert _{H\left( A\right) }=\left\Vert Au\right\Vert
+\left\Vert u\right\Vert . 
\]

The symmetric operator $A$ is positive defined if $D\left( A\right) $ is
dense in $H$ and there exists a positive constant $C_{A}$ depend only on $A$
such that 
\[
\left( Au,u\right) \geq C_{A}\left\Vert u\right\Vert ^{2}. 
\]

It is known that see e.g. $\left[ \text{24, \S 1.15.1}\right] $) there exist
fractional powers\ $A^{\theta }$ of the positive defined operator $A.$

Let $\left[ A,B\right] $ be a commutator operator, i.e. 
\[
\left[ A,B\right] =AB-BA 
\]%
for linear operators $A$ and $B.$

Sometimes we use one and the same symbol $C$ without distinction in order to
denote positive constants which may differ from each other even in a single
context. When we want to specify the dependence of such a constant on a
parameter, say $\alpha $, we write $C_{\alpha }$.

\begin{center}
\textbf{2}.\textbf{1}. \textbf{Main results for} \textbf{absract Scr\"{o}%
dinger equation}
\end{center}

Consider the problem $\left( 1.1\right) $. \ Here, 
\[
X=L^{2}\left( R^{n};H\right) \text{, }X\left( A\right) =L^{2}\left(
R^{n};H\left( A\right) \right) ,\text{\ }Y^{k}=W^{2,k}\left( R^{n};H\right) 
\text{, }k\in \mathbb{N}. 
\]

\textbf{Definition} \textbf{2.1}. A function $u\in L^{\infty }\left(
0,T;H\left( A\right) \right) $ is called a local weak solution to $\left(
1.1\right) $ on $\left( 0,T\right) $ if $u$ belongs to $L^{\infty }\left(
0,T;H\left( A\right) \right) \cap W^{1,2}\left( 0,T;H\right) $ and satisfies 
$\left( 1.1\right) $ in the sense of $L^{\infty }\left( 0,T;H\left( A\right)
\right) .$ In particular, if $\left( 0,T\right) $\ coincides with $\mathbb{R}
$, then $u$ is called a global weak solution to $\left( 1.1\right) .$

If the solution of $\left( 1.1\right) $ belongs to $C\left( \left[ 0,T\right]
;H\left( A\right) \right) \cap W^{2,2}\left( 0,T;H\right) ,$ then its called
a stronge solution.

\bigskip \textbf{Condition 1. }Assume: (1)\textbf{\ }$A=A\left( x\right) $
is a symmetric operator in Hilbert space $H$ with independent on $x$ domain $%
D\left( A\right) $ that is dense on $H;$

(2) $\frac{\partial A}{\partial x_{k}}$ are symmetric operators in Hilbert
space $H$ with independent on $x$ domain $D\left( \frac{\partial A}{\partial
x_{k}}\right) =D\left( A\right) $ and 
\[
\dsum\limits_{k=1}^{n}\left( x_{k}\left[ A\frac{\partial f}{\partial x_{k}}-%
\frac{\partial A}{\partial x_{k}}f\right] ,f\right) _{X}\geq 0,\text{ } 
\]%
for each $f\in L^{\infty }\left( 0,T;Y^{1}\left( A\right) \right) $;

(3) 
\[
\dsum\limits_{k=1}^{n}\left( x_{k}\left[ A\frac{\partial f}{\partial x_{k}}+%
\frac{\partial A}{\partial x_{k}}f\right] ,f\right) _{X}\geq 0\text{ } 
\]%
for each $f\in L^{\infty }\left( 0,T;Y^{1}\left( A\right) \right) $;

(4) $iA$ generates a Schr\"{o}dinger grop and $V\left( x,t\right) \in
L\left( H\right) $ for $\left( x,t\right) \in R^{n}\times \left[ 0,1\right]
; $

(5) there is a constant $C_{0}>0$ so that 
\[
\func{Im}\left( \left( A+V\right) \upsilon ,\upsilon \right) _{H}\geq
C_{0}\mu \left( x,t\right) \left\Vert \upsilon \left( x,t\right) \right\Vert
^{2}, 
\]%
for $x\in R^{n},$ $t\in \left[ 0,T\right] ,$ $T\in \left( 0.\right. 1\left.
{}\right] $ and $\upsilon \in D\left( A\right) $, where $\mu $ is a positive
function in $L^{1}\left( 0,T;L^{\infty }\left( R^{n}\right) \right) ;$

(6)\ 
\[
\left\Vert V\right\Vert _{L^{\infty }\left( R^{n}\times \left( 0,1\right)
;L\left( H\right) \right) }\leq C\text{, }\lim\limits_{R\rightarrow \infty
}\left\Vert V\right\Vert _{L_{t}^{1}L_{x}^{\infty }\left( L\left( H\right)
\right) }=0, 
\]%
where 
\[
\text{ }L_{t}^{1}L_{x}^{\infty }\left( L\left( H\right) \right) =L^{1}\left(
0,1;L^{\infty }\left( R^{n}/O_{r}\right) ;L\left( H\right) \right) . 
\]

Let $A=A\left( x\right) $ is a symmetric operator in Hilbert space $H$ with
independent on $x$ domain $D\left( A\right) $ and%
\[
X\left( A\right) =L^{2}\left( R^{n};H\left( A\right) \right) ,\text{ }%
Y^{k}\left( A\right) =W^{2,k}\left( R^{n};H\left( A\right) \right) \text{, }%
k\in \mathbb{N} 
\]%
\ Our main result in this paper is the following:

\textbf{Theorem 1. }Assume the Condition 1 holds and there exist constants $%
a_{0},$ $a_{1},$ $a_{2}>0$ such that for any $k\in \mathbb{Z}^{+}$ a
solution $u\in C\left( \left[ 0,1\right] ;X\left( A\right) \right) $ of $%
\left( 1.1\right) $ satisfy 
\begin{equation}
\dint\limits_{R^{n}}\left\Vert u\left( x,0\right) \right\Vert
^{2}e^{2a_{0}\left\vert x\right\vert ^{p}}dx<\infty ,\text{ for }p\in \left(
1,2\right) ,  \tag{2.1}
\end{equation}%
\begin{equation}
\dint\limits_{R^{n}}\left\Vert u\left( x,1\right) \right\Vert
^{2}e^{2k\left\vert x\right\vert ^{p}}dx<a_{2}e^{2a_{1}k^{\frac{q}{q-p}}},%
\text{ }\frac{1}{p}+\frac{1}{q}=1.  \tag{2.2}
\end{equation}%
Moreover, there exists $M_{p}>0$ such that

\begin{equation}
a_{0}a_{1}^{p-2}>M_{p}.  \tag{2.3}
\end{equation}%
Then $u\left( x,t\right) \equiv 0.$

\bigskip \textbf{Corollary 1. }Assume the Condition1 holds and 
\[
\lim\limits_{\left\vert R\right\vert \rightarrow \infty
}\dint\limits_{0}^{1}\sup\limits_{\left\vert x\right\vert >R}\left\Vert
V\left( x,t\right) \right\Vert _{L\left( H\right) }dt=0. 
\]%
There exist positive constants $\alpha $, $\beta $ such that a solution $%
u\in C\left( \left[ 0,1\right] ;X\left( A\right) \right) $ of $\left(
1.1\right) $ satisfy 
\begin{equation}
\dint\limits_{R^{n}}\left\Vert u\left( x,0\right) \right\Vert
^{2}e^{2\left\vert \alpha x\right\vert
^{p}/p}dx+\dint\limits_{R^{n}}\left\Vert u\left( x,1\right) \right\Vert
^{2}e^{\frac{2\left\vert \beta x\right\vert ^{q}}{q}}dx<\infty ,\text{ } 
\tag{2.4}
\end{equation}%
with 
\[
p\in \left( 1,2\right) \text{, }\frac{1}{p}+\frac{1}{q}=1 
\]%
and there exists $N_{p}>0$ such that 
\begin{equation}
\alpha \beta >N_{p}.\text{ }  \tag{2.5}
\end{equation}%
Then $u\left( x,t\right) \equiv 0.$

As a direct consequence of Corollary1 we have the following result regarding
the uniqueness of solutions for nonlinear equation $\left( 1.2\right) $:

\bigskip \textbf{Theorem 2. }Assume the Condition 1 holds and $u_{1},$ $%
u_{2}\in C\left( \left[ 0,1\right] ;Y^{2,k}\left( A\right) \right) $ strong
solutions of $(1.2)$ with $k\in \mathbb{Z}^{+},$ $k>\frac{n}{2}.$ Suppose $%
F:H\times H\rightarrow H,$ $F\in C^{k}$, $F\left( 0\right) =\partial
_{u}F\left( 0\right) =\partial _{\bar{u}}F\left( 0\right) =0$ and there
exist positive constants $\alpha $, $\beta $ such that%
\begin{equation}
e^{\frac{\left\vert \alpha x\right\vert ^{p}}{p}}\left( u_{1}\left(
.,0\right) -u_{2}\left( .,0\right) \right) \in X,\text{ }e^{^{\frac{%
\left\vert \beta x\right\vert ^{q}}{q}}}\left( u_{1}\left( .,0\right)
-u_{2}\left( .,0\right) \right) \in X,  \tag{2.6}
\end{equation}%
with 
\[
p\in \left( 1,2\right) \text{, }\frac{1}{p}+\frac{1}{q}=1 
\]%
and there exists $N_{p}>0$ such that 
\begin{equation}
\alpha \beta >N_{p}.\text{ }  \tag{2.7}
\end{equation}

Then $u_{1}\left( x,t\right) \equiv u_{2}\left( x,t\right) .$

\textbf{Corollary 2.} Assume the Condition 1 hold and there exist positive
constants $\alpha $ and $\beta $ such that a solution $u\in C\left( \left[
0,1\right] ;X\left( A\right) \right) $ of $\left( 1.1\right) $ satisfy 
\begin{equation}
\dint\limits_{R^{n}}\left\Vert u\left( x,0\right) \right\Vert ^{2}e^{\frac{%
2\left\vert \alpha x_{j}\right\vert ^{p}}{p}}dx+\dint\limits_{R^{n}}\left%
\Vert u\left( x,1\right) \right\Vert ^{2}e^{\frac{2\left\vert \beta
x_{j}\right\vert ^{q}}{q}}dx<\infty ,\text{ }  \tag{2.8}
\end{equation}%
for $j=1,2,...,n$ and $p\in \left( 1,2\right) $, $\frac{1}{p}+\frac{1}{q}=1.$
Moreover, there exists $N_{p}>0$ such that 
\begin{equation}
\alpha \beta >N_{p}.\text{ }  \tag{2.9}
\end{equation}%
Then $u\left( x,t\right) \equiv 0.$

\textbf{Remark 2.1. }The Theorem 2 still holds, with different constant $%
N_{p}>0$, if one replaces the hypothesis $(2.6)$ by

\[
e^{\left\vert \alpha x_{j}\right\vert ^{p}/p}\left( u_{1}\left( .,0\right)
-u_{2}\left( .,0\right) \right) \in X\left( A\right) ,\text{ }e^{\left\vert
\beta x_{j}\right\vert ^{q}/q}\left( u_{1}\left( .,0\right) -u_{2}\left(
.,0\right) \right) \in X\left( A\right) , 
\]%
for $j=1,2,...,n.$

Next, we shall extend the method used in the proof Theorem 3 to study the
blow the nonlinear Schr\"{o}dinger equations

\begin{equation}
i\partial _{t}u+\Delta u+Au+F\left( u,\bar{u}\right) u=0,\text{ }x\in R^{n},%
\text{ }t\in \left[ 0,1\right] ,  \tag{2.10}
\end{equation}%
\bigskip where $A$ is a linear operator in a Hilbert space $H.$

Let $u\left( x,t\right) $ be a solution of the equation $\left( 2.10\right)
. $ Then it can be shown that the function%
\begin{equation}
u\left( x,t\right) =U\left( x,1-t\right) u\left( \frac{x}{1-t},\frac{t}{1-t}%
\right) \text{, }  \tag{2.11}
\end{equation}%
is a solution of the focussing $L^{2}$-critical solution of abstract
Schredinger equation%
\begin{equation}
i\partial _{t}u+\Delta u+Au+\left\Vert u\right\Vert ^{\frac{4}{n}}u=0,\text{ 
}x\in R^{n},\text{ }t\in \left[ 0,1\right]  \tag{2.12}
\end{equation}%
which blows up at time $t=1,$ where $U\left( x,t\right) $ is a fundamental
solution of the Schr\"{o}dinger equation%
\[
i\partial _{t}u+\Delta u+Au=0,\text{ }x\in R^{n},\text{ }t\in \left[ 0,1%
\right] , 
\]%
i.e. 
\[
U\left( x,t\right) =t^{-n/2}\exp \left\{ i\left( A+\left\vert x\right\vert
^{2}\right) /4t\right\} . 
\]%
By using the above result we wil prove the following main result

\textbf{Theorem 3. }Assume the Condition 1 holds and there exist positive
constants $b_{0}$ and $\theta $ such that a solution $u\in C\left( \left[
-1,1\right] ;X\left( A\right) \right) $ of $\left( 2.10\right) $ satisfied $%
\left\Vert F\left( u,\bar{u}\right) \right\Vert \leq b_{0}\left\Vert
u\right\Vert ^{\theta }$ for $\left\Vert u\right\Vert >1$. Suppose 
\[
\left\Vert u\left( .,t\right) \right\Vert _{X}=\left\Vert u\left( .,0\right)
\right\Vert _{X}=\left\Vert u_{0}\right\Vert _{X}=a,\text{ }t\in \left(
-1,1\right) 
\]%
and 
\begin{equation}
\left\Vert u\left( x,t\right) \right\Vert \leq \left( 1-t\right) ^{-\frac{n}{%
2}}Q\left( \frac{\left\vert x\right\vert }{1-t}\right) ,\text{ }t\in \left(
-1,1\right)  \tag{2.14}
\end{equation}%
\ where, 
\begin{equation}
Q\left( x\right) =b_{1}^{-\frac{n}{2}}e^{-b_{2}\left\vert x\right\vert ^{p}}%
\text{, }b_{1}\text{, }b_{2}>0\text{, }p>1;  \tag{2.15}
\end{equation}%
If $p>p\left( \theta \right) =\frac{2\left( \theta n-2\right) }{\left(
\theta n-1\right) },$ then $a\equiv 0.$

\begin{center}
\textbf{2.2. Some auxiliary results}
\end{center}

First of all, we generalize the result G. W. Morgan (see e.g $\left[ 7\right]
$) about Morgan type uncertainty principle for Fourier transform.

\textbf{Lemma 2.0. }Let \textbf{\ }$f\left( x\right) \in L^{1}\left(
R^{n};H\right) \cap X$ and%
\[
\dint\limits_{R^{n}}\dint\limits_{R^{n}}\left\Vert f\left( x\right)
\right\Vert \left\Vert \hat{f}\left( \xi \right) \right\Vert e^{\left\vert
x.\zeta \right\vert }\text{ }dxd\xi <\infty \text{.} 
\]

Then $f\left( x\right) \equiv 0.$

In particular, using Young's inequality this implies:

\textbf{Result 2.1.} Let%
\[
f\left( x\right) \in L^{1}\left( R^{n};H\right) \cap X,\text{ }p\in \left(
1,2\right) ,\text{ }\frac{1}{p}+\frac{1}{q}=1,\alpha ,\beta >0 
\]%
and 
\[
\dint\limits_{R^{n}}\left\Vert f\left( x\right) \right\Vert e^{\frac{\alpha
^{p}\left\vert x\right\vert ^{p}}{p}}dx+\dint\limits_{R^{n}}\left\Vert \hat{f%
}\left( \xi \right) \right\Vert e^{\frac{\beta ^{q}\left\vert \xi
\right\vert ^{q}}{q}}d\xi <\infty \text{, }\alpha \beta >1. 
\]%
Then $f\left( x\right) \equiv 0.$

The Morgan type uncertainty principle, in terms of the solution of the free
Schredinger equation will be as:

Let%
\[
u_{0}\left( .\right) \in L^{1}\left( R^{n};H\right) \cap X\left( A\right) 
\]%
and 
\[
\dint\limits_{R^{n}}\left\Vert u_{0}\left( x\right) \right\Vert e^{\frac{%
\alpha ^{p}\left\vert x\right\vert ^{p}}{p}}dx+\dint\limits_{R^{n}}\left%
\Vert e^{it\left( \Delta +A\right) }u_{0}\left( x\right) \right\Vert e^{%
\frac{\beta ^{q}\left\vert \xi \right\vert }{q\left( 2t\right) ^{q}}%
^{q}}dx<\infty \text{, }\alpha \beta >1 
\]%
for some $t\neq 0$. Then $u_{0}\left( x\right) \equiv 0.$

Let $L^{2}\left( H\right) =L^{2}\left( \left( 0,1\right) \times
R^{n};H\right) .$

Consider the abstract Schredinger equation%
\begin{equation}
\partial _{t}u=\left( a+ib\right) \left[ \Delta u+A\left( x\right) u+V\left(
x,t\right) u+F\left( x,t\right) \right] ,\text{ }x\in R^{n},\text{ }t\in %
\left[ 0,1\right] ,  \tag{2.16}
\end{equation}%
where $a$, $b$ are real numbers, $A$ is a linear operator, $V\left(
x,t\right) $ is a given potential operator function in $H$ and $F\left(
x,t\right) $ is a given $H$-valued function.

Let\ $S=S\left( t\right) $ be a symmetric, $K=K\left( t\right) $ be a
skew-symmetric operators in $H$, \ $f=f\left( .,t\right) \in X$ for $t\in %
\left[ 0,1\right] $ and 
\[
\text{ }Q\left( t\right) =\left( f,f\right) _{X}\text{, }D\left( t\right)
=\left( Sf,f\right) _{X},\text{ }N\left( t\right) =D\left( t\right)
Q^{-1}\left( t\right) , 
\]%
\[
\text{ }\partial _{t}S=S_{t}\text{ and }\left\vert \nabla \upsilon
\right\vert _{H}^{2}=\dsum\limits_{k=1}^{n}\left\Vert \frac{\partial
\upsilon }{\partial x_{k}}\right\Vert _{H}^{2}\text{ for }\upsilon \in
W^{1,2}\left( R^{n};H\right) . 
\]

Let 
\[
\eta \left( x\right) =\frac{\gamma a\left\vert x\right\vert ^{2}}{a+4\gamma
\left( a^{2}+b^{2}\right) T},\text{ }L^{1,\infty }\left( 0,T;H\right)
=L^{1}\left( 0,T;L^{\infty }\left( R^{n};H\right) \right) . 
\]

\textbf{Lemma 2.1. \ }Let the Condition 1 holds$.$ Assume that $u\in
L^{\infty }\left( \left[ 0,1\right] ;X\left( A\right) \right) \cap
L^{2}\left( 0,1;Y^{1}\right) $ is a solution of $\left( 2.16\right) $ for $%
a\in \mathbb{R}_{+}$ and $b\in \mathbb{R}.$ Then, 
\[
e^{-M_{T}}\left\Vert e^{\eta \left( x\right) }u\left( .,T\right) \right\Vert
_{X}\leq 
\]%
\[
\left\Vert e^{\gamma \left\vert x\right\vert ^{2}}u\left( .,0\right)
\right\Vert _{X}+\left( a^{2}+b^{2}\right) \left\Vert e^{\eta \left(
x\right) }F\left( .,t\right) \right\Vert _{L^{1}\left( 0,T;X\right)
}+\left\Vert A^{\frac{1}{2}}u\left( .,t\right) \right\Vert _{X}^{2}, 
\]%
uniformly in $t\in \left[ 0,1\right] ,$\ when $\gamma \geq 0,$ $0\leq T\leq
1 $ and%
\[
M_{T}=\left\Vert a\left( \func{Re}V\right) ^{+}-b\func{Im}V\right\Vert
_{L^{1,\infty }\left( 0,T;H\right) }. 
\]

\textbf{Proof. }Set $\upsilon $ $=e^{\varphi }u$, where $\varphi $ is a
real-valued function to be chosen later. The function $\upsilon =\upsilon
\left( x,t\right) $ verifies 
\begin{equation}
\partial _{t}\upsilon =S\upsilon +K\upsilon +\left( a+ib\right) \left(
Vf+e^{\gamma \varphi }F\right) \text{ in }R^{n}\times \left[ 0,1\right] 
\text{,}  \tag{2.17}
\end{equation}%
where $S,$ $K$ are symmetric and skew-symmetric operator, respectively given
by

\[
S=aA_{1}-ib\gamma B_{1}+\varphi _{t}+a\func{Re}V-b\func{Im}V,\text{ }%
K=ibA_{1}-a\gamma B_{1}+ 
\]%
\begin{equation}
i\left( b\func{Re}\upsilon +a\func{Im}\upsilon \right) ,  \tag{2.18}
\end{equation}%
here 
\[
A_{1}=\Delta +A\left( x\right) +\gamma ^{2}\left\vert \nabla \varphi
\right\vert ^{2},\text{ }B_{1}=2\nabla \varphi .\nabla +\Delta \varphi . 
\]%
Formally, \ 
\[
\partial _{t}\left\Vert \upsilon \right\Vert _{X}^{2}=2\func{Re}\left(
S\upsilon ,\upsilon \right) _{X}+2\func{Re}\left( \left( a+ib\right)
e^{\gamma \varphi }F,\upsilon \right) _{X} 
\]%
for $t\geq 0.$ Again a formal integration by parts gives that 
\[
\func{Re}\left( S\upsilon ,\upsilon \right)
_{X}=-a\dint\limits_{R^{n}}\left\vert \nabla \upsilon \right\vert
_{H}^{2}dx+a\dint\limits_{R^{n}}\left( A\upsilon ,\upsilon \right)
dx+a\gamma ^{2}\dint\limits_{R^{n}}\left\vert \nabla \varphi \right\vert
^{2}\left\Vert \upsilon \left( x,t\right) \right\Vert ^{2}dx+ 
\]%
\[
\dint\limits_{R^{n}}\varphi _{t}\left\Vert \upsilon \left( x,t\right)
\right\Vert ^{2}dx+b\func{Im}\dint\limits_{R^{n}}\left( \bar{\upsilon}\nabla
\varphi ,\nabla \upsilon \right) dx+\dint\limits_{R^{n}}\left( \left( a\func{%
Re}V-b\func{Im}V\right) \upsilon ,\upsilon \right) dx. 
\]%
By Cauchy-Schwarz's inequality and in view of assumption on $A,$ $V$ we get 
\[
\partial _{t}\left\Vert \upsilon \right\Vert _{X}^{2}\leq 2a\left\Vert A^{%
\frac{1}{2}}\upsilon \right\Vert _{X}^{2}+2\left\Vert \left( a\func{Re}V-b%
\func{Im}V\right) \upsilon \right\Vert _{L^{\infty }\left( L\left( H\right)
\right) }\left\Vert \upsilon \left( .,t\right) \right\Vert _{X}^{2}+ 
\]%
\begin{equation}
2\sqrt{a^{2}+b^{2}}\left\Vert e^{\gamma \varphi }F\left( .,t\right)
\right\Vert _{X}\left\Vert \upsilon \left( .,t\right) \right\Vert _{X}, 
\tag{2.19}
\end{equation}%
when 
\begin{equation}
\left( a+\frac{b^{2}}{a}\right) \left\vert \nabla \upsilon \right\vert
^{2}+\varphi _{t}\leq 0,\text{ for }x\in R^{n},\text{ }t\geq 0.  \tag{2.20}
\end{equation}

\bigskip Let 
\begin{equation}
\phi _{r}\left( x\right) =\left\{ 
\begin{array}{c}
\left\vert x\right\vert ^{2},\text{ \ \ }\left\vert x\right\vert \leq r \\ 
r^{2},\text{ \ }\left\vert x\right\vert >r%
\end{array}%
\right. \text{, }\varphi _{\rho ,r}=d\left( t\right) \theta _{\rho }\ast
\phi _{r}\left( x\right) \text{, }\upsilon _{\rho ,r}=e^{\varphi _{\rho
,r}}u,  \tag{2.21}
\end{equation}%
where $\theta _{\rho }$ is a radial mollifier and 
\[
d\left( t\right) =\frac{\gamma a}{a+4\gamma \left( a^{2}+b^{2}\right) t}. 
\]%
Then by reasoning as in the end of proof $\left[ \text{8, Lemma 1}\right] ,$
from $\left( 2.19\right) $ we obtain that the estimate holds 
\[
\left\Vert \upsilon \left( .,T\right) \right\Vert _{X}\leq e^{M_{T}}\left[
\left\Vert e^{\gamma \left\vert x\right\vert ^{2}}u\left( .,T\right)
\right\Vert _{X}\right. +2a\left\Vert A^{\frac{1}{2}}u\left( .,t\right)
\right\Vert _{X}^{2}+ 
\]%
\[
\left. 2\sqrt{a^{2}+b^{2}}\left\Vert e^{\varphi _{\rho ,r}}F\right\Vert
_{L^{1}\left( 0,T;X\right) }\right] 
\]%
holds uniformly in $\rho $, $r$ and $t\in \left[ 0,1\right] .$ Lemma 2.1
follows after letting $\rho $ tend to zero and $r$ to infinity.

\textbf{Lemma 2.2. \ }Let $S=S\left( x,t\right) $ be a symmetric operator, $%
K=K\left( x,t\right) $ be skew-symmetric in Hilbert space $H$ with
independent on $x$ domain $D\left( S\right) $ and $D\left( K\right) ,$
respectively$.$\textbf{\ }Assume\ $G\left( .,t\right) \in L^{2}\left(
R^{n}\right) $ is a positive funtion and $f(x,t)$ is a $H-$valued reasonable
function. Then, 
\[
\frac{d^{2}}{dt^{2}}Q\left( t\right) =2\partial _{t}\func{Re}\left( \partial
_{t}f-Sf-Kf,f\right) _{X}+2\left( S_{t}f+\left[ S,K\right] f,f\right) _{X}+ 
\]

\[
\left\Vert \partial _{t}f-Sf+Kf\right\Vert _{X}^{2}-\left\Vert \partial
_{t}f-Sf-Kf\right\Vert _{X}^{2} 
\]%
and 
\[
\partial _{t}N\left( t\right) \geq Q^{-1}\left( t\right) \left[ \left(
S_{t}f+\left[ S,K\right] f,f\right) _{X}-\frac{1}{2}\left\Vert \partial
_{t}f-Sf-Kf\right\Vert _{X}^{2}\right] . 
\]%
Moreover, if 
\[
\left\Vert \partial _{t}f-Sf-Kf\right\Vert _{H}\leq M_{1}\left\Vert
f\right\Vert _{H}+G\left( x,t\right) ,\text{ }S_{t}+\left[ S,K\right] \geq
-M_{0}\text{ } 
\]%
for $x\in R^{N}$, $t\in \left[ 0,1\right] $ and%
\[
M_{2}=\sup\limits_{t\in \left[ 0,1\right] }\left\Vert G\left( .,t\right)
\right\Vert _{L^{2}\left( R^{n}\right) }\left\Vert f\left( .,t\right)
\right\Vert _{X}^{-1}<\infty . 
\]

Then $Q\left( t\right) $ is logarithmically convex in $[0,1]$ and there is a
constant $M$ such that 
\[
Q\left( t\right) \leq e^{M\left(
M_{0}+M_{1}+M_{2}+M_{1}^{2}+M_{2}^{2}\right) }Q^{1-t}\left( 0\right)
Q^{t}\left( 1\right) \text{, }0\leq t\leq 1. 
\]

\textbf{Proof.} It is clear that 
\[
\dot{Q}\left( t\right) =2\func{Re}\left( \frac{d}{dt}f,f\right) _{X}=2\func{%
Re}\left( \frac{d}{dt}f-Sf-Kf,f\right) _{X}+2\left( Sf,f\right) _{X}. 
\]

Also, 
\[
\dot{Q}\left( t\right) =\func{Re}\left( \frac{d}{dt}f+Sf,f\right) _{X}+\func{%
Re}\left( \frac{d}{dt}f-Sf,f\right) _{X}, 
\]

\[
\frac{d}{dt}D\left( t\right) =\frac{1}{2}\func{Re}\left( \frac{d}{dt}%
f+Sf,f\right) _{X}-\frac{1}{2}\func{Re}\left( \frac{d}{dt}f-Sf,f\right)
_{X}. 
\]

\bigskip Then by reasoning as in $\left[ \text{8, Lemma 2}\right] $ we
obtain the assertion.

Consider the following abstract Schredinger equation 
\begin{equation}
\partial _{t}u=\left( a+ib\right) \left[ \Delta u+A\left( x\right) u+V\left(
x,t\right) u\right] +F\left( x,t\right) ,\text{ }x\in R^{n},\text{ }t\in %
\left[ 0,1\right] ,  \tag{2.22}
\end{equation}

Let 
\[
\Phi \left( A,V\right) \upsilon =a\func{Re}\left( \left( A+V\right) \upsilon
,\upsilon \right) -b\func{Im}\left( \left( A+V\right) \upsilon ,\upsilon
\right) , 
\]%
\[
\text{ for }\upsilon =\upsilon \left( .,t\right) \in D\left( A\right) . 
\]

\bigskip \textbf{Lemma 2.3. \ }Let the Condition 1 holds. Assume that $a$, $%
\gamma >0,$ $b\in \mathbb{R}.$ Let $u\in L^{\infty }\left( 0,1;X\left(
A\right) \right) \cap L^{2}\left( 0,1;Y^{1}\right) $ be solution of the
equation $\left( 2.22\right) $ and 
\[
\left\vert \Phi \left( A,V\right) u\left( x,t\right) \right\vert \leq
C_{0}\eta \left( x,t\right) \left\Vert u\left( x,t\right) \right\Vert ^{2} 
\]%
for $x\in R^{n},$ $t\in \left[ 0,1\right] $ and $u\in D\left( A\right) $,
where $\eta $ is a positive function in $L^{1}\left( 0,T;L^{\infty }\left(
R^{n}\right) \right) $. Moreover, suppose%
\[
\sup\limits_{t\in \left[ 0,1\right] }\left\Vert V\left( .,t\right)
\right\Vert _{B}\leq M_{1}\text{, }\left\Vert e^{\gamma \left\vert
x\right\vert ^{2}}u\left( .,0\right) \right\Vert _{X}<\infty ,\text{ }%
\left\Vert e^{\gamma \left\vert x\right\vert ^{2}}u\left( .,1\right)
\right\Vert _{X}<\infty \text{ } 
\]%
and 
\[
M_{2}=\sup\limits_{t\in \left[ 0,1\right] }\frac{\left\Vert e^{\gamma
\left\vert x\right\vert ^{2}}F\left( .,t\right) \right\Vert _{X}}{\left\Vert
u\right\Vert _{X}}<\infty . 
\]

Then, $e^{\gamma \left\vert x\right\vert ^{2}}u\left( .,t\right) $ is
logarithmically convex in $[0,1]$ and there is a constant $N$ such that%
\begin{equation}
\left\Vert e^{\gamma \left\vert x\right\vert ^{2}}u\left( .,t\right)
\right\Vert _{X}\leq e^{NM\left( a,b\right) }\left\Vert e^{\gamma \left\vert
x\right\vert ^{2}}u\left( .,0\right) \right\Vert _{X}^{1-t}\left\Vert
e^{\gamma \left\vert x\right\vert ^{2}}u\left( .,1\right) \right\Vert
_{X}^{t}  \tag{2.23}
\end{equation}%
where 
\[
M\left( a,b\right) =\left( a^{2}+b^{2}\right) \left( \gamma
M_{1}^{2}+M_{2}^{2}\right) +\sqrt{a^{2}+b^{2}}\left( M_{1}+M_{2}\right) 
\]%
when $0\leq t\leq 1.$

\textbf{Proof. }Let $f=e^{\gamma \varphi }u,$ where $\varphi $ is a
real-valued function to be chosen. The function $f\left( x,t\right) $
verifies%
\begin{equation}
\partial _{t}f=Sf+Kf+\left( a+ib\right) \left( Vf+e^{\gamma \varphi
}F\right) \text{ in }R^{n}\times \left[ 0,1\right] \text{,}  \tag{2.24}
\end{equation}%
where $S$, $K$ are symmetric and skew-symmetric operator, respectively given
by%
\[
S=aA_{1}-ib\gamma B_{1}+\gamma \varphi _{t},\text{ }K=ibA_{1}-a\gamma B_{1}, 
\]%
where 
\[
A_{1}=\Delta +A+\gamma ^{2}\left\vert \nabla \varphi \right\vert ^{2},\text{ 
}B_{1}=2\nabla \varphi .\nabla +\Delta \varphi . 
\]%
A calculation shows that 
\begin{equation}
S_{t}+\left[ S,K\right] =\gamma \partial _{t}^{2}\varphi +2\gamma
^{2}a\nabla \varphi .\nabla \varphi _{t}-2ib\gamma \left( 2\nabla \varphi
_{t}.\nabla +\Delta \varphi _{t}\right) -  \tag{2.25}
\end{equation}%
\[
\gamma \left( a^{2}+b^{2}\right) \left[ 4\nabla .\left( D^{2}\varphi \nabla
\right) -4\gamma ^{2}D^{2}\varphi \nabla \varphi +\Delta ^{2}\varphi \right]
+2\left[ A\left( x\right) \nabla \varphi .\nabla -\nabla \varphi .\nabla A%
\right] . 
\]%
If we put $\varphi =\left\vert x\right\vert ^{2}$, then $\left( 2.25\right) $
reduce the following%
\[
S_{t}+\left[ S,K\right] =-\gamma \left( a^{2}+b^{2}\right) \left[ 8\Delta
-32\gamma ^{2}\left\vert x\right\vert ^{2}\right] +2\left[ A\left( x\right)
\nabla \varphi .\nabla -\nabla \varphi .\nabla A\right] . 
\]%
Moreover by assumtion (2) of Condition 1, 
\[
\left( S_{t}f+\left[ S,K\right] f,f\right) =\gamma \left( a^{2}+b^{2}\right)
\dint\limits_{R^{n}}\left( 8\left\Vert \left\vert \nabla f\right\vert
\right\Vert ^{2}+32\gamma ^{2}\left\vert x\right\vert ^{2}\left\Vert
f\right\Vert ^{2}\right) dx+ 
\]%
\[
4\dsum\limits_{k=1}^{n}\dint\limits_{R^{n}}x_{k}\left( \left[ A\frac{%
\partial f}{\partial x_{k}}-\frac{\partial A}{\partial x_{k}}f\right]
,f\right) dx\geq 0\text{ for }f\in L^{\infty }\left( 0,T;Y^{1}\left(
A\right) \right) . 
\]%
This identity, the condition on $A=A\left( x\right) ,$ $V=V\left( x,t\right) 
$ and $\left( 2.24\right) $ imply that%
\begin{equation}
\text{ }S_{t}+\left[ S,K\right] \geq 0,  \tag{2.26}
\end{equation}%
\[
\left\Vert \partial _{t}f-Sf-Kf\right\Vert _{X}\leq \sqrt{a^{2}+b^{2}}\left(
M_{1}\left\Vert f\right\Vert _{X}+e^{\gamma \varphi }\left\Vert F\right\Vert
_{X}\right) \text{.} 
\]%
If we knew that the quantities and calculations involved in the proof of
Lemma 2.2. ( this fact are derived \i n a similar way as in $\left[ \text{8,
Lemma 2}\right] $) were finite and correct, when $f=e^{\gamma \left\vert
x\right\vert ^{2}}u$, we would have the logarithmic convexity of $%
Q(t)=\left\Vert e^{\gamma \left\vert x\right\vert ^{2}}u\left( .,t\right)
\right\Vert _{X}$ and get $(2.23)$ from Lemma 2.2. To justify the validity
of the previous arguments now we consider the function $\varphi _{a}\left(
x\right) $ constructed in $\left[ \text{8, Lemma 3}\right] $, i.e., 
\[
\varphi _{a}\left( x\right) =\left\{ 
\begin{array}{c}
\left\vert x\right\vert ^{2},\text{ \ \ }\left\vert x\right\vert \leq 1 \\ 
\left( 2-a\right) ^{-1}\left( 2\left\vert x\right\vert ^{2-a}-a\right) ,%
\text{ \ }\left\vert x\right\vert \geq 1%
\end{array}%
\right. 
\]%
and replace $\varphi =\left\vert x\right\vert ^{2}$ by $\varphi _{a,\rho
}=\theta _{\rho }\ast \varphi _{a},$ where $a,$ $\rho \in \left( 0,1\right) $
and $\theta \in C_{0}^{\infty }\left( R^{n}\right) $ is a radial function.
From the proof $\left[ \text{8, Lemma 3}\right] $ we get that $\varphi
_{a,\rho }$ is convex function and 
\begin{equation}
\left\Vert \Delta \varphi _{a,\rho }\right\Vert _{\infty }\leq C\left(
n,\rho \right) a.  \tag{2.27}
\end{equation}%
Put then, $f_{a,\rho }=e^{\gamma \varphi _{a,\rho }}u$ for $u\in X$ and $%
Q_{a,\rho }\left( t\right) =$ $\left\Vert f_{a,\rho }\left( .,t\right)
\right\Vert _{X}^{2}$ in Lemma 2.2. The decay bound in Lemma 2.1 and the
interior regularity for solutions of $(2.23)$ can now be used qualitatively
to make sure that the quantities or calculations involved in the proof of $%
\left[ \text{8, Lemma 2}\right] $ ( with replicing $L^{2}$ by $X)$ are
finite and correct for $f_{a,\rho }$. In this case, $f_{a,\rho }$ verifies%
\[
\partial _{t}f_{a,\rho }=S^{a,\rho }f_{a,\rho }+K^{a,\rho }f_{a,\rho
}+\left( a+ib\right) \left[ A\left( x\right) f_{a,\rho }+e^{\gamma \varphi
_{a,\rho }}F\left( x,t\right) \right] 
\]%
with symmetric and skew-symmetric operators $S^{a,\rho }$ and $K^{a,\rho }$
given by $(2.24)$ with $\varphi $ replaced by $\varphi _{a,\rho }$. The
formula $\left( 2\text{.}25\right) $, the convexity of $\varphi _{a,\rho }$,
the bounds $(2.26)$ and $(2.27)$ imply that the inequalities%
\[
\text{ }S_{t}^{a,\rho }+\left[ S^{a,\rho },K^{a,\rho }\right] \geq 0, 
\]%
\[
\left\Vert \partial _{t}f^{a,\rho }-Sf^{a,\rho }-K^{a,\rho }f\right\Vert
_{X}\leq \sqrt{a^{2}+b^{2}}\left( M_{1}\left\Vert f^{a,\rho }\right\Vert
_{X}+e^{\gamma \varphi _{a,\rho }}\left\Vert F\right\Vert _{X}\right) \text{.%
} 
\]%
hold and $M_{2}\left( a,\rho \right) \leq e^{C\left( n\right) \rho
^{2}}M_{2},$ when $a,\rho \in \left( 0,1\right) .$ Particularly, $Q_{a,\rho
} $ is logarithmically convex in $[0,1]$ and%
\begin{equation}
Q_{a,\rho }\left( t\right) \leq e^{N\left[ \left( a^{2}+b^{2}\right) \left(
M_{1}^{2}+M_{2}^{2}\right) +\left( M_{1}+M_{2}\right) \sqrt{a^{2}+b^{2}}%
\right] }Q_{a,\rho }^{1-t}\left( 0\right) Q_{a,\rho }^{t}\left( 1\right) . 
\tag{2.28}
\end{equation}

Then, $(2.23)$ follows after taking first the limit, when a tends to zero in 
$(2.28)$ and then, when $\rho $ tends to zero.

\begin{center}
\textbf{3. Some properties of solutions of the abstract Schredinger equations%
}
\end{center}

Let

\begin{center}
\[
\sigma \left( t\right) =\left[ \alpha \left( 1-t\right) +\beta t\right]
^{-1},\text{ }\eta \left( x,t\right) =\left( \alpha -\beta \right) )|x|^{2}%
\left[ 4i(\alpha (1-t)+\beta t)\right] ^{-1},\text{ } 
\]
\end{center}

\[
\nu \left( s\right) =\left[ \gamma \alpha \beta \sigma ^{2}\left( s\right) +%
\frac{\left( \alpha -\beta \right) a}{4\left( a^{2}+b^{2}\right) }\sigma
\left( s\right) \right] ,\text{ }\phi \left( x,t\right) =\text{ }\frac{%
\gamma a\left\vert x\right\vert ^{2}}{a+4\gamma \left( a^{2}+b^{2}\right) t}%
. 
\]

Let $u=u\left( x,s\right) $ be a solution of the equation 
\[
\partial _{s}u=i\left[ \Delta u+Au+V\left( y,s\right) u+F\left( y,s\right) %
\right] ,\text{ }y\in R^{n},\text{ }s\in \left[ 0,1\right] . 
\]%
and $a+ib\neq 0$, $\gamma \in \mathbb{R}$, $\alpha $, $\beta \in \mathbb{R}%
_{+}$. Set 
\begin{equation}
\tilde{u}\left( x,t\right) =\left( \sqrt{\alpha \beta }\sigma \left(
t\right) \right) ^{\frac{n}{2}}u\left( \sqrt{\alpha \beta }x\sigma \left(
t\right) ,\beta t\sigma \left( t\right) \right) e^{\eta }.  \tag{3.1}
\end{equation}

\bigskip Then, $\tilde{u}\left( x,t\right) $\ verifies the equation 
\begin{equation}
\partial _{t}\tilde{u}=i\left[ \Delta \tilde{u}+A\left( x\right) \tilde{u}+%
\tilde{V}\left( x,t\right) \tilde{u}+\tilde{F}\left( x,t\right) \right] ,%
\text{ }x\in R^{n},\text{ }t\in \left[ 0,1\right]  \tag{3.2}
\end{equation}%
with 
\begin{equation}
\tilde{V}\left( x,t\right) =\alpha \beta \sigma ^{2}\left( t\right) V\left( 
\sqrt{\alpha \beta }x\sigma \left( t\right) ,\beta t\sigma \left( t\right)
\right) ,  \tag{3.3}
\end{equation}

\begin{equation}
\text{ }\tilde{F}\left( x,t\right) =\left( \sqrt{\alpha \beta }\sigma \left(
t\right) \right) ^{\frac{n}{2}+2}\left( \sqrt{\alpha \beta }x\sigma \left(
t\right) ,\beta t\sigma \left( t\right) \right) .  \tag{3.4}
\end{equation}%
Moreover, 
\begin{equation}
\left\Vert e^{\gamma \left\vert x\right\vert ^{2}}\tilde{F}\left( .,t\right)
\right\Vert _{X}=\alpha \beta \sigma ^{2}\left( t\right) e^{\nu \left\vert
y\right\vert ^{2}}\left\Vert F\left( s\right) \right\Vert _{X}\text{, }%
\left\Vert e^{\gamma \left\vert x\right\vert ^{2}}\tilde{u}\left( .,t\right)
\right\Vert _{X}=e^{\nu \left\vert y\right\vert ^{2}}\left\Vert u\left(
s\right) \right\Vert _{X}  \tag{3.5}
\end{equation}%
when $s=\beta t\sigma \left( t\right) $.

\textbf{Remark 3.1. }Let $\beta =\beta \left( k\right) .$ By assumption we
have 
\[
\left\Vert e^{a_{0}\left\vert x\right\vert ^{p}}u\left( x,0\right)
\right\Vert _{X}=a_{0}, 
\]%
\begin{equation}
\left\Vert e^{k\left\vert x\right\vert ^{p}}u\left( x,0\right) \right\Vert
_{X}=a_{k}\leq a_{2}e^{2a_{1}k^{\frac{q}{q-p}}}=a_{2}e^{2a_{1}k^{\frac{1}{2-p%
}}}.  \tag{3.6}
\end{equation}%
Thus, for $\gamma =\gamma (k)$ $\in \lbrack 0,\infty )$ to be chosen later,
one has%
\begin{equation}
\left\Vert e^{\gamma \left\vert x\right\vert ^{p}}\tilde{u}_{k}\left(
x,0\right) \right\Vert _{X}=\left\Vert e^{\gamma \left( \frac{\alpha }{\beta 
}\right) ^{p/2}\left\vert x\right\vert ^{p}}u_{k}\left( x,0\right)
\right\Vert _{X}=b_{0},  \tag{3.7}
\end{equation}

\[
\left\Vert e^{\gamma \left\vert x\right\vert ^{p}}\tilde{u}_{k}\left(
x,1\right) \right\Vert _{X}=\left\Vert e^{\gamma \left( \frac{\beta }{\alpha 
}\right) ^{p/2}\left\vert x\right\vert ^{p}}u_{k}\left( x,1\right)
\right\Vert _{X}=a_{k}. 
\]%
Let we take 
\[
\gamma \left( \frac{\alpha }{\beta }\right) ^{p/2}=a_{0}\text{ and }\gamma
\left( \frac{\beta }{\alpha }\right) ^{p/2}=k, 
\]%
i.e. 
\begin{equation}
\gamma =\left( ka_{0}\right) ^{\frac{1}{2}}\text{, }\beta =k^{\frac{1}{p}},%
\text{ }\alpha =a_{0}^{\frac{1}{p}}.  \tag{3.8}
\end{equation}

Let 
\[
M=\dint\limits_{0}^{1}\left\Vert V\left( .,t\right) \right\Vert _{L^{\infty
}\left( R^{n};H\right) }dt=\dint\limits_{0}^{1}\left\Vert V\left( .,s\right)
\right\Vert _{L^{\infty }\left( R^{n};H\right) }ds. 
\]

From $\left( 3.2\right) $, using energy estimates it follows 
\begin{equation}
e^{-M}\left\Vert u\left( .,0\right) \right\Vert _{X}\leq \left\Vert u\left(
.,t\right) \right\Vert _{X}=\left\Vert \tilde{u}\left( .,s\right)
\right\Vert _{X}\leq e^{M}\left\Vert u\left( .,0\right) \right\Vert _{X}%
\text{, }t,s\in \left[ 0,1\right] ,  \tag{3.9}
\end{equation}%
where 
\[
s=\beta t\sigma \left( t\right) . 
\]

Consider the following problem%
\begin{equation}
i\partial _{t}u+\Delta u+A\left( x\right) u+V\left( x,t\right) u+F\left(
x,t\right) =0,\text{ }x\in R^{n},\text{ }t\in \left[ 0,1\right] ,  \tag{3.10}
\end{equation}

\[
u\left( x,0\right) =u_{0}\left( x\right) ,\text{ } 
\]%
where $A=A\left( x\right) $ is a linear operator$,$ $V\left( x,t\right) $ is
a given potential operator function in a Hilbert space $H$ and $F$ is a $H$%
-valued function.

Let as define operator valued integral operators in $L^{p}\left( \Omega
;H\right) $. Let $k$: $R^{n}\backslash \left\{ 0\right\} \rightarrow L\left(
H\right) .$ We say $k\left( x\right) $ is a $L\left( H\right) $-valued
Calderon-Zygmund kernel ($C-Z$ kernel) if $k\in C^{\infty }\left(
R^{n}\backslash \left\{ 0\right\} ,L\left( H\right) \right) ,$ $k$ is
homogenous of degree $-n,$ $\dint\limits_{B}k\left( x\right) d\sigma =0,$
where%
\[
B=\left\{ x\in R^{n}\text{: }\left\vert x\right\vert =1\right\} . 
\]%
For $f\in L^{p}\left( \Omega ;H\right) ,$ $p\in \left( 1,\infty \right) ,$ $%
a\in L^{\infty }\left( R^{n}\right) ,$ and $x\in \Omega $ we set the
Calderon-Zygmund operator 
\[
K_{\varepsilon }f=\dint\limits_{\left\vert x-y\right\vert >\varepsilon ,y\in
\Omega }k\left( x,y\right) f\left( y\right) dy,\text{ }Kf=\lim\limits_{%
\varepsilon \rightarrow 0}K_{\varepsilon }f 
\]%
and commutator operator%
\[
\left[ K;a\right] f=a\left( x\right) Kf\left( x\right) -K\left( af\right)
\left( x\right) = 
\]%
\[
\lim\limits_{\varepsilon \rightarrow 0}\dint\limits_{\left\vert
x-y\right\vert >\varepsilon ,y\in \Omega }k\left( x,y\right) \left[ a\left(
x\right) -a\left( y\right) \right] f\left( y\right) dy. 
\]

By using Calder\'{o}n's first commutator estimates $\left[ 4\right] ,$
convolution operators on abstract functions $\left[ 2\right] $ and abstract
commutator theorem in $\left[ 20\right] $ we obtain the following result:

\textbf{Theorem A}$_{2}.$ Assume $k(.)$ is $L\left( H\right) -$valued $C-Z$
kernel that have locally integrable first-order derivatives in $\left\vert
x\right\vert >0$, and 
\[
\left\Vert k\left( x,y\right) -k\left( x^{\prime },y\right) \right\Vert
_{L\left( E\right) }\leq M\left\vert x-x^{\prime }\right\vert \left\vert
x-y\right\vert ^{-\left( n+1\right) }\text{ for }\left\vert x-y\right\vert
>2\left\vert x-x^{\prime }\right\vert . 
\]%
Let $a(.)$ have first-order derivatives in $L^{r}\left( R^{n}\right) $,$%
1<r\leq \infty $. Then for $p$, $q\in \left( 1,\infty \right) ,$ $%
q^{-1}=p^{-1}+r^{-1}$ the following estimates hold%
\[
\left\Vert \left[ K;a\right] \partial _{x_{j}}f\right\Vert _{L^{q}\left(
R^{n};H\right) }\leq C\left\Vert f\right\Vert _{L^{p}\left( R^{n};H\right)
}, 
\]%
\[
\left\Vert \partial _{x_{j}}\left[ K;a\right] f\right\Vert _{L^{q}\left(
R^{n};H\right) }\leq C\left\Vert f\right\Vert _{L^{p}\left( R^{n};H\right)
}, 
\]%
for $f\in C_{0}^{\infty }\left( R^{n};H\right) ,$\ where the constant $C>0$
is independent of $f.$

Let $X_{\gamma }$ denote the weithed Lebesque space $L_{\gamma }^{2}\left(
R^{n}\text{:}H\right) $ with $\gamma \left( x\right) =e^{2\lambda .x}.$ By
following $\left[ \text{5, Lemma 2.1}\right] $ let us show:

\textbf{Lemma 3.1.} Assume that the Condition 1 holds and there exists $%
\varepsilon _{0}>0$ such that%
\[
\left\Vert V\right\Vert _{L_{t}^{1}L_{x}^{\infty }\left( R^{n}\times \left[
0,1\right] ;L\left( H\right) \right) }<\varepsilon _{0}. 
\]%
Moreover, suppose $u\in C\left( \left[ 0,1\right] ;X\left( A\right) \right) $
is a stronge solution of $\left( 3.10\right) $ with 
\[
u_{0},\text{ }u_{1}=u\left( x,1\right) \in X_{\gamma },\text{ }F\in
L^{1}\left( 0,1;X_{\gamma }\right) 
\]%
for $\gamma \left( x\right) =e^{2\lambda .x}$ $\ $and for some $\lambda \in
R^{n}.$ Then there exists a positive constant $M_{0}=M_{0}\left(
n,A,H\right) $ independent of $\lambda $ such that%
\begin{equation}
\sup\limits_{t\in \left[ 0,1\right] }\left\Vert u\left( .,t\right)
\right\Vert _{X_{\gamma }}\leq M_{0}\left[ \left\Vert u_{0}\right\Vert
_{X_{\gamma }}+\left\Vert u_{1}\right\Vert _{X_{\gamma
}}+\dint\limits_{0}^{1}\left\Vert F\left( .,t\right) \right\Vert _{X_{\gamma
}}dt\right] .  \tag{3.11}
\end{equation}

\textbf{Proof. }First, we consider the case, when $\gamma \left( x\right)
=\beta \left( x\right) =e^{2\beta x_{1}}.$ Without loss of generality we
shall assume $\beta >0.$ Let $\varphi _{n}\in C^{\infty }\left( \mathbb{R}%
\right) $ such that $\varphi _{n}\left( \tau \right) =1$, $\tau \leq n$ and $%
\varphi _{n}\left( \tau \right) =0$ for $\tau \geq 10n$ with $0\leq \varphi
_{n}\leq 1,$ $\left\vert \varphi _{n}^{\left( j\right) }\left( \tau \right)
\right\vert \leq C_{j}\tau n^{-j}.$ Let 
\[
\theta _{n}(\tau )=\beta \dint\limits_{0}^{\tau }\varphi _{n}^{2}\left(
s\right) ds 
\]%
so that $\theta _{n}\in C^{\infty }\left( \mathbb{R}\right) $ nondecreasing
with $\theta _{n}(\tau )=\beta \tau $ for $\tau <n,$ $\theta _{n}(\tau
)=C_{n}\beta $ for $\tau >10n$ and 
\begin{equation}
\theta _{n}^{\prime }(\tau )=\beta \varphi _{n}^{2}\left( \tau \right) \leq
\beta ,\text{ }\theta _{n}^{\left( j\right) }(\tau )=\beta C_{j}n^{1-j},%
\text{ }j=1,2,....  \tag{3.12}
\end{equation}

\bigskip Let $\phi _{n}\left( \tau \right) =\exp \left( 2\theta _{n}(\tau
)\right) $ so that $\phi _{n}\left( \tau \right) \leq \exp \left( 2\beta
\tau \right) $ and $\phi _{n}\left( \tau \right) \rightarrow \exp \left(
2\beta \tau \right) $ for $n\rightarrow \infty .$ Let $u\left( x,t\right) $
be a solution of the equation $\left( 3.10\right) $, then one gets the
equation $\upsilon _{n}\left( x,t\right) =\phi _{n}\left( x_{1}\right)
u\left( x,t\right) $ satisfies the following 
\begin{equation}
i\partial _{t}\upsilon _{n}+\Delta \upsilon _{n}+A\upsilon _{n}=V_{n}\left(
x,t\right) \upsilon _{n}+\phi _{n}\left( x_{1}\right) F\left( x,t\right) , 
\tag{3.13}
\end{equation}%
where 
\[
V_{n}\left( x,t\right) \upsilon _{n}=V\left( x,t\right) \upsilon _{n}+4\beta
\varphi _{n}\left( x_{1}\right) \partial _{x_{1}}\upsilon _{n}+\left[ 4\beta
\varphi _{n}\left( x_{1}\right) \varphi _{n}^{\prime }\left( x_{1}\right)
-4\beta ^{2}\upsilon _{n}^{4}\right] \upsilon _{n}. 
\]%
Now, we consider a new function

\[
w_{n}\left( x,t\right) =e^{\mu }\upsilon _{n}\left( x,t\right) ,\text{ }\mu
=-i4\beta ^{2}\varphi _{n}^{4}\left( x_{1}\right) t. 
\]%
Then from $\left( 3.13\right) $ we get%
\[
i\partial _{t}w_{n}+\Delta w_{n}+Aw_{n}=\tilde{V}_{n}\left( x,t\right) w_{n}+%
\tilde{F}_{n}\left( x,t\right) , 
\]%
where%
\[
\tilde{V}_{n}\left( x,t\right) w_{n}=V\left( x,t\right) w_{n}+h\left(
x_{1},t\right) +a^{2}\left( x_{1}\right) \partial _{x_{1}}w_{n}+itb\left(
x_{1}\right) \partial _{x_{1}}w_{n} 
\]%
\[
\tilde{F}_{n}\left( x,t\right) =e^{\mu }\phi _{n}\left( x_{1}\right) F\left(
x,t\right) , 
\]%
when 
\[
h\left( x_{1},t\right) =\left( i16\beta ^{2}\varphi _{n}^{3}\varphi
_{n}^{\prime }t\right) ^{2}+i48\beta ^{2}\varphi _{n}^{2}\left( \varphi
_{n}^{\prime }\right) ^{2}t+i16\beta ^{2}\varphi _{n}^{3}\varphi
_{n}^{\left( 2\right) }t+ 
\]

\[
4\beta \varphi _{n}\varphi _{n}^{\prime }+i64\beta ^{2}\varphi
_{n}^{3}\varphi _{n}^{\prime }t,\text{ }a^{2}=4\beta \varphi _{n}^{2}\left(
x_{1}\right) \text{, }b=-32\beta ^{2}\varphi _{n}^{3}\varphi _{n}^{\prime }. 
\]

It is clear to see that 
\[
\left\Vert \partial _{x_{1}}^{j}h\left( x_{1},t\right) \right\Vert
_{L^{\infty }\left( \mathbb{R}\times \left[ 0,1\right] \right) }\leq
C_{j}n^{-\left( j+1\right) },\text{ }j=1,2,..., 
\]

\begin{equation}
a^{2}\left( x_{1}\right) \geq 0,\text{ }\left\Vert \partial
_{x_{1}}^{j}a\left( x_{1}\right) \right\Vert _{L^{\infty }\left( \mathbb{R}%
\right) }\leq C_{j}n^{-j},\text{ }j=1,2,...,  \tag{3.15}
\end{equation}

\[
\left\Vert \partial _{x_{1}}^{j}b\left( x_{1}\right) \right\Vert _{L^{\infty
}\left( \mathbb{R}\right) }\leq C_{j}n^{-j},\text{ }j=1,2,.... 
\]

\bigskip Then by reasoning as in $\left[ \text{5, Lemma 2.1}\right] $ and by
using the properties of symmetric operators $A$ and $V$, we obtain 
\[
\partial _{t}\left\vert \left( P_{\varepsilon }P_{+}w_{n},\upsilon \right)
\right\vert ^{2}+2\func{Im}\left( \Delta P_{\varepsilon }P_{+}w_{n},\upsilon
\right) \left( \bar{K}\left( P_{\varepsilon }P_{+}w_{n},\upsilon \right)
\right) + 
\]

\[
2\func{Im}\left( AP_{\varepsilon }P_{+}w_{n},\upsilon \right) \left( \bar{K}%
\left( P_{\varepsilon }P_{+}w_{n},\upsilon \right) \right) =2\func{Im}\left(
P_{\varepsilon }P_{+}\left( Vw_{n}\right) ,\upsilon \right) \left( \bar{K}%
\left( P_{\varepsilon }P_{+}w_{n},\upsilon \right) \right) + 
\]

\[
2\func{Im}\left( \left( P_{\varepsilon }P_{+}hw_{n}\right) ,\upsilon \right)
\left( \bar{K}\left( P_{\varepsilon }P_{+}w_{n},\upsilon \right) \right) +2%
\func{Im}\left( \left( P_{\varepsilon }P_{+}w_{n},\upsilon \right) \right) 
\bar{K}\left( P_{\varepsilon }P_{+}\left( \tilde{F}_{n}\right) ,\upsilon
\right) + 
\]

\begin{equation}
2\func{Im}\left( \left( P_{\varepsilon }P_{+}a^{2}\left( x_{1}\right)
\partial _{x_{1}}w_{n}\right) ,\upsilon \right) \left( \bar{K}\left(
P_{\varepsilon }P_{+}w_{n},\upsilon \right) \right) +  \tag{3.16}
\end{equation}%
\[
2\func{Im}\left( \left( P_{\varepsilon }P_{+}ib\left( x_{1}\right) \partial
_{x_{1}}w_{n}\right) ,\upsilon \right) \left( \bar{K}\left( P_{\varepsilon
}P_{+}w_{n},\upsilon \right) \right) . 
\]

\bigskip Since for all $n\in \mathbb{Z}^{+}$, $w_{n}\left( .\right) \in X$, $%
\tilde{F}_{n}\left( .,t\right) \in X$ and for a.e. $t\in \left[ 0,1\right] $
by integrating both sides of $\left( 3.16\right) $ on $R^{n}$ we get 
\[
\dint\limits_{R^{n}}\func{Im}\left( \Delta P_{\varepsilon
}P_{+}w_{n},\upsilon \right) \left( \bar{K}\left( P_{\varepsilon
}P_{+}w_{n},\upsilon \right) \right) dx=0. 
\]

It is clear to see that%
\[
\left( u,\upsilon \right) _{X}=\left( u\left( .\right) ,\upsilon \left(
.\right) _{H}\right) _{L^{2}\left( R^{n}\right) }\text{, for }u\text{, }%
\upsilon \in X. 
\]

\bigskip Then applying the Cauchy-Schwartz and Holder inequalites for a.e. $%
t\in \left[ 0,1\right] $ we obtain 
\begin{equation}
\dint\limits_{R^{n}}\func{Im}\left( P_{\varepsilon }P_{+}Vw_{n},\upsilon
\right) \left( \bar{K}\left( P_{\varepsilon }P_{+}w_{n},\upsilon \right)
\right) dx\leq C\left\Vert V\right\Vert _{B}\left\Vert w_{n}\right\Vert
_{X}^{2}\left\Vert \upsilon \right\Vert _{X}^{2},  \tag{3.17}
\end{equation}%
\begin{equation}
\dint\limits_{R^{n}}\func{Im}\left( P_{\varepsilon }P_{+}hw_{n},\upsilon
\right) \bar{K}\left( P_{\varepsilon }P_{+}w_{n},\upsilon \right) dx\leq
C\left\Vert h\right\Vert _{L^{\infty }}\left\Vert w_{n}\right\Vert
_{X}^{2}\left\Vert \upsilon \right\Vert _{X}^{2},  \tag{3.18}
\end{equation}%
\begin{equation}
\dint\limits_{R^{n}}\func{Im}\left( P_{\varepsilon }P_{+}\tilde{F}%
_{n}w_{n},\upsilon \right) \left( \bar{K}\left( P_{\varepsilon
}P_{+}w_{n},\upsilon \right) \right) dx\leq C\left\Vert \tilde{F}%
_{n}\right\Vert _{X}\left\Vert w_{n}\right\Vert _{X}^{2}\left\Vert \upsilon
\right\Vert _{X}^{2}.  \tag{3.19}
\end{equation}

Moreover, again applying the Cauchy-Schwartz and Holder inequalities due to
symmetricity of the operator $A,$ for a.e. $t\in \left[ 0,1\right] $ we get

\begin{equation}
\dint\limits_{R^{n}}\func{Im}\left( AP_{\varepsilon }P_{+}w_{n},\upsilon
\right) \left( \bar{K}\left( P_{\varepsilon }P_{+}w_{n},\upsilon \right)
\right) dx\leq C\left\Vert Aw_{n}\right\Vert _{X}\left\Vert w_{n}\right\Vert
_{X}\left\Vert \upsilon \right\Vert _{X}^{2}  \tag{3.20}
\end{equation}%
where, the constant $C$ in $\left( 3.18\right) -\left( 3.20\right) $ is
independent of $\upsilon \in C_{0}^{\infty }\left( R^{n};H\right) ,$ $%
\varepsilon \in \left( 0,\left. 1\right] \right. $ and $n\in \mathbb{Z}^{+}.$
Since $C_{0}^{\infty }\left( R^{n};H\right) $ is dense in $X,$ from $\left(
3.17\right) $-$\left( 3.20\right) $ in view of operator theory in Hilbert
spaces, we obtain the following%
\[
\left\vert \dint\limits_{R^{n}}\func{Im}\left( P_{\varepsilon }P_{+}Vw_{n},%
\bar{K}P_{\varepsilon }P_{+}w_{n}\right) dx\right\vert \leq C\left\Vert
V\right\Vert _{B}\left\Vert w_{n}\right\Vert _{X}^{2}, 
\]%
\[
\left\vert \dint\limits_{R^{n}}\func{Im}\left( P_{\varepsilon }P_{+}hw_{n},%
\bar{K}P_{\varepsilon }P_{+}w_{n}\right) dx\right\vert \leq C\left\Vert
h\right\Vert _{L^{\infty }}\left\Vert w_{n}\right\Vert _{X}^{2}\leq C\frac{1%
}{n}\left\Vert w_{n}\right\Vert _{X}^{2}, 
\]%
\begin{equation}
\left\vert \dint\limits_{R^{n}}\func{Im}\left( P_{\varepsilon }P_{+}\tilde{F}%
_{n}w_{n},\bar{K}P_{\varepsilon }P_{+}w_{n}\right) dx\right\vert \leq
C\left\Vert \tilde{F}_{n}\right\Vert _{X}\left\Vert w_{n}\right\Vert
_{X}^{2},  \tag{3.21}
\end{equation}

\[
\left\vert \dint\limits_{R^{n}}\func{Im}\left( AP_{\varepsilon }P_{+}w_{n},%
\bar{K}P_{\varepsilon }P_{+}w_{n}\right) dx\right\vert \leq C\left\Vert
Aw_{n}\right\Vert _{X}\left\Vert w_{n}\right\Vert _{X}^{2} 
\]

For bounding the last two terms in $\left( 3.16\right) $ we will use the
abstract version of Calder\'{o}n's first commutator estimates $\left[ 4%
\right] .$ Really, by Cauchy-Schvartz inequality and in view of Theorem A$%
_{2}$ we get 
\begin{equation}
\left\Vert \left( \left[ P_{\pm };a\right] \partial _{x_{1}}f,\upsilon
\right) \right\Vert _{X}\leq C\left\Vert \partial _{x_{1}}a\right\Vert
_{L^{\infty }}\left\Vert f\right\Vert _{X}\left\Vert \upsilon \right\Vert
_{X},\text{ }  \tag{3.22}
\end{equation}%
\begin{equation}
\left\Vert \partial _{x_{1}}\left( \left[ P_{\pm };a\right] f,\upsilon
\right) \right\Vert _{X}\leq C\left\Vert \partial _{x_{1}}a\right\Vert
_{L^{\infty }}\left\Vert f\right\Vert _{X}\left\Vert \upsilon \right\Vert
_{X},\text{ }  \tag{3.23}
\end{equation}

\bigskip Also, from the calculus of pseudodifferential operators with
operator coefficients (see e.g. $\left[ 5\right] $ ) and the inequality $%
(3.15)$, we have

\begin{equation}
\left\Vert \left( \left[ P_{\varepsilon };a\right] \partial
_{x_{1}}f,\upsilon \right) \right\Vert _{X}\leq \frac{C}{n}\left\Vert
f\right\Vert _{X}\left\Vert \upsilon \right\Vert _{X},\text{ }  \tag{3.24}
\end{equation}%
\begin{equation}
\left\Vert \partial _{x_{1}}\left( \left[ P_{\varepsilon };a\right]
f,\upsilon \right) \right\Vert _{X}\leq \frac{C}{n}\left\Vert f\right\Vert
_{X}\left\Vert \upsilon \right\Vert _{X},\text{ }  \tag{3.25}
\end{equation}%
where the constant $C$ in $\left( 3.22\right) -(3.25)$ is independent of $%
\varepsilon \in \lbrack 0,1]$ and $n.$ We remark that estimates $%
(3.22)-(3.25)$ also hold with $b\left( x_{1}\right) $ replacing $a(x_{1})$.
Since $C_{0}^{\infty }\left( R^{n};H\right) $ is dense in $X,$ from $\left(
3.22\right) $-$\left( 3.25\right) $ in view of operator theory in Hilbert
spaces, we obtain

\[
\left\Vert \left[ P_{\pm };a\right] \partial _{x_{1}}f\right\Vert _{X}\leq
C\left\Vert \partial _{x_{1}}a\right\Vert _{L^{\infty }}\left\Vert
f\right\Vert _{X},\text{ }\left\Vert \partial _{x_{1}}\left[ P_{\pm };a%
\right] f\right\Vert _{X}\leq C\left\Vert \partial _{x_{1}}a\right\Vert
_{L^{\infty }}\left\Vert f\right\Vert _{X}, 
\]%
\begin{equation}
\left\Vert \left( \left[ P_{\varepsilon };a\right] \partial
_{x_{1}}f,\upsilon \right) \right\Vert _{X}\leq \frac{C}{n}\left\Vert
f\right\Vert _{X},\text{ }\left\Vert \partial _{x_{1}}\left( \left[
P_{\varepsilon };a\right] f,\upsilon \right) \right\Vert _{X}\leq \frac{C}{n}%
\left\Vert f\right\Vert _{X}\text{ ,}  \tag{3.26}
\end{equation}%
and the same estimates $\left( 3.26\right) $ with $b(x_{1}$) replacing $%
a(x_{1})$. By reasoning as in $\left[ \text{6, Lemma 2.1}\right] $ (claim 1
and 2 ) from $\left( 3.26\right) $\ we obtain 
\begin{equation}
\left\vert \func{Im}\left( \left( P_{\varepsilon }P_{+}a^{2}\left(
x_{1}\right) \partial _{x_{1}}w_{n},\bar{K}P_{\varepsilon }P_{+}w_{n}\right)
\right) \right\vert \leq O\left( n^{-1}\left\Vert w_{n}\right\Vert
_{X}\right) ,  \tag{3.27}
\end{equation}

\[
\left\vert \func{Im}\left( \left( P_{\varepsilon }P_{+}b\left( x_{1}\right)
\partial _{x_{1}}w_{n},\bar{K}P_{\varepsilon }P_{+}w_{n}\right) \right)
\right\vert \leq O\left( n^{-1}\left\Vert w_{n}\right\Vert _{X}\right) . 
\]

Now, the estimates $\left( 3.21\right) $ and $\left( 3.27\right) $ implay
the assertion.

\begin{center}
\textbf{4.} \textbf{Proof of Theorem 1}
\end{center}

\textbf{\ }We will apply Lemma 3.1 to a solution of the equation $(3.2)$.
Since $0<\alpha <\beta =$ $\beta (k)$ for $k>k_{0}$ it follows that $\alpha
\leq \sigma \left( t\right) \leq \beta $ for any $t\in \lbrack 0,1]$.
Therefore if $y=\sqrt{\alpha \beta }x\sigma \left( t\right) $, then from $%
\left( 3.8\right) $ we get

\begin{equation}
\sqrt{\alpha \beta ^{-1}}\left\vert x\right\vert \leq \left\vert
y\right\vert \sqrt{\alpha ^{-1}\beta }\left\vert x\right\vert =\left(
ka_{0}^{-1}\right) ^{\frac{1}{2p}}\left\vert x\right\vert  \tag{4.1}
\end{equation}

\bigskip Thus, 
\begin{equation}
\left\Vert \alpha \beta \sigma ^{2}\left( t\right) V\left( \sqrt{\alpha
\beta }x\sigma \left( t\right) ,\beta t\sigma \left( t\right) \right)
\right\Vert _{L\left( H\right) }\leq \alpha ^{-1}\beta \left\Vert
V\right\Vert _{B}=\left( ka_{0}^{-1}\right) ^{\frac{1}{p}}\left\Vert
V\right\Vert _{B}  \tag{4.2}
\end{equation}%
and so, 
\begin{equation}
\left\Vert \tilde{V}\left( .,t\right) \right\Vert _{L^{\infty }\left(
R^{n};H\right) }\leq \left( ka_{0}^{-1}\right) ^{\frac{1}{p}}\left\Vert
V\left( .,t\right) \right\Vert _{L^{\infty }\left( R^{n};H\right) }. 
\tag{4.3}
\end{equation}

Also, for $s=\beta t\sigma \left( t\right) $ it is clear that 
\begin{equation}
\frac{ds}{dt}=\alpha \beta \sigma ^{2}\left( t\right) \text{, }dt=\left(
\alpha \beta \right) ^{-1}\sigma ^{-2}\left( t\right) ds.  \tag{4.4}
\end{equation}%
Therefore, 
\[
\dint\limits_{0}^{1}\left\Vert \tilde{V}\left( .,t\right) \right\Vert
_{L^{\infty }\left( R^{n};H\right) }dt=\dint\limits_{0}^{1}\left\Vert
V\left( .,s\right) \right\Vert _{L^{\infty }\left( R^{n};H\right) }ds, 
\]%
and from $\left( 4.1\right) $ we get 
\begin{equation}
\dint\limits_{0}^{1}\left\Vert \tilde{V}\left( .,t\right) \right\Vert
_{L^{\infty }\left( \left\vert x\right\vert >r;H\right)
}dt=\dint\limits_{0}^{1}\left\Vert V\left( .,s\right) \right\Vert
_{L^{\infty }\left( \left\vert y\right\vert >\varkappa ;H\right) }ds, 
\tag{4.5}
\end{equation}%
where 
\[
\varkappa =\left( a_{0}k^{-1}\right) ^{\frac{1}{2p}}r. 
\]%
So, if $\dint\limits_{0}^{1}\left\Vert V\left( .,s\right) \right\Vert
_{L^{\infty }\left( \left\vert y\right\vert >\varkappa ;H\right)
}ds<\varepsilon _{0}$ then, 
\[
\dint\limits_{0}^{1}\left\Vert \tilde{V}\left( .,t\right) \right\Vert
_{L^{\infty }\left( \left\vert y\right\vert >r;H\right) }ds<\varepsilon _{0},%
\text{ for }r=\varkappa \left( ka_{0}^{-1}\right) ^{\frac{1}{2p}} 
\]%
and we can applay Lemma 3.1 to the equation $\left( 3.2\right) $ with%
\[
\tilde{V}\mathbb{=}\tilde{V}_{\chi \left( \left\vert x\right\vert >r\right)
}\left( x,t\right) \text{, }\tilde{F}\mathbb{=}\tilde{V}_{\chi \left(
\left\vert x\right\vert <R\right) }\left( x,t\right) \tilde{u}\left(
x,t\right) 
\]%
to get the following estimate 
\[
\sup\limits_{t\in \left[ 0,1\right] }\left\Vert e^{\nu }\tilde{u}\left(
.,t\right) \right\Vert _{X}\leq M_{0}\left( \left\Vert e^{\nu }\tilde{u}%
\left( .,0\right) \right\Vert _{X}+\left\Vert e^{\nu }\tilde{u}\left(
.,1\right) \right\Vert _{X}\right) + 
\]%
\begin{equation}
M_{0}e^{M}e^{\nu _{0}}\left\Vert \tilde{V}\right\Vert _{B}\left\Vert u\left(
.,0\right) \right\Vert _{X},  \tag{4.6}
\end{equation}%
where $M$ a positive constant defined in Remark 2.1 and 
\[
B=L^{\infty }\left( R^{n}\times \left[ 0,1\right] ;L\left( H\right) \right) 
\text{, }\nu =\left( 2p\right) ^{\frac{1}{p}}\gamma ^{\frac{1}{p}}\lambda .%
\frac{x}{2},\text{ }\nu _{0}=\left\vert \lambda \right\vert \left( 2p\right)
^{\frac{1}{p}}\gamma ^{\frac{1}{p}}\frac{r}{2}. 
\]%
From $\left( 4.6\right) $ we have 
\[
\sup\limits_{t\in \left[ 0,1\right] }\dint\limits_{R^{n}}\left\Vert e^{\nu }%
\tilde{u}\left( .,t\right) \right\Vert _{H}^{2}dx\leq
M_{0}\dint\limits_{R^{n}}e^{\nu }\left( \left\Vert \tilde{u}\left(
.,0\right) \right\Vert _{H}^{2}+\left\Vert \tilde{u}\left( .,1\right)
\right\Vert _{H}^{2}\right) dx+ 
\]%
\[
M_{0}e^{M}e^{\left\vert \lambda \right\vert \left( 2p\right) ^{\frac{1}{p}%
}\gamma ^{\frac{1}{p}r}}\left\Vert \tilde{V}\right\Vert _{B}\left\Vert
u\left( .,0\right) \right\Vert _{X}^{2}, 
\]%
and multiply the above inequality by $e^{\left\vert \lambda \right\vert
/q}\left\vert \lambda \right\vert ^{n\left( q-2\right) /2}$, integrate in $%
\lambda $ and in $x$, use Fubini theorem and the following formula 
\begin{equation}
e^{\gamma \left\vert x\right\vert ^{p}/p}\thickapprox
\dint\limits_{R^{n}}e^{\gamma ^{\frac{1}{p}}\lambda .x-\left\vert \lambda
\right\vert ^{q/q}}\left\vert \lambda \right\vert ^{n\left( q-2\right)
/2}d\lambda ,  \tag{4.7}
\end{equation}%
proven in $\left[ \text{7, Appendx}\right] $ to obtain 
\[
\dint\limits_{\left\vert x\right\vert >1}e^{2\gamma \left\vert x\right\vert
^{p}}\left\Vert \tilde{u}\left( .,t\right) \right\Vert _{H}^{2}dx\leq
M_{0}\dint\limits_{R^{n}}e^{2\gamma \left\vert x\right\vert ^{p}}\left(
\left\Vert \tilde{u}\left( .,0\right) \right\Vert _{H}^{2}+\left\Vert \tilde{%
u}\left( .,1\right) \right\Vert _{H}^{2}\right) dx+ 
\]%
\begin{equation}
M_{0}e^{2M}e^{2\gamma r^{p}}r^{C_{p}}\left\Vert \tilde{V}\right\Vert
_{B}\left\Vert u\left( .,0\right) \right\Vert _{X}^{2}.  \tag{4.8}
\end{equation}%
Hence, the esimates $\left( 3.6\right) $, $\left( 3.8\right) $, $\left(
3.9\right) ,$ $\left( 4.3\right) $ and $\left( 4.8\right) $ imply 
\[
\sup\limits_{t\in \left[ 0,1\right] }\left\Vert e^{\gamma \left\vert
x\right\vert ^{p}}\tilde{u}\left( .,t\right) \right\Vert _{X}\leq
M_{0}\left( \left\Vert e^{\gamma \left\vert x\right\vert ^{p}}\tilde{u}%
\left( .,0\right) \right\Vert _{X}+\left\Vert e^{\gamma \left\vert
x\right\vert ^{p}}\tilde{u}\left( .,1\right) \right\Vert _{X}\right) + 
\]%
\begin{equation}
M_{0}e^{M}e^{\gamma }\left\Vert u\left( .,0\right) \right\Vert
_{X}+M_{0}\left( ka_{0}^{-1}\right) ^{C_{p}}e^{M}e^{\varkappa ^{p}\gamma ^{%
\frac{1}{2}}\left( ka_{0}^{-1}\right) }\left\Vert u\left( .,0\right)
\right\Vert _{X}\left\Vert V\right\Vert _{B}\leq  \tag{4.9}
\end{equation}%
\[
M_{0}\left( a_{0}+a_{k}\right) +M_{0}e^{M}\left\Vert u\left( .,0\right)
\right\Vert _{X}\left( e^{\gamma }+\left( ka_{0}^{-1}\right)
^{C_{p}}\left\Vert V\right\Vert _{B}\right) e^{k\varkappa ^{p}}\leq 
\]%
\[
M_{0}a_{k}=M_{0}e^{a_{1}k^{1/2-p}}\text{ for }k>k_{0}\left( M_{0}\right) 
\text{ sufficiently large.} 
\]

Next, we shall obtain bounds for the $\nabla \tilde{u}$. Let $\tilde{\gamma}=%
\frac{\gamma }{2}$ and $\varphi $ be a strictly convex complex valued
function on compact sets of $R^{n}$, radial such that (see $[7]$)%
\[
D^{2}\varphi \geq p(p-1)|x|^{(p-2)},\text{ for }|x|\geq 1, 
\]%
\[
\varphi \geq 0,\text{ }\left\Vert \partial ^{\alpha }\varphi \right\Vert
_{L^{\infty }}\leq C\text{, }2\leq |\alpha |\leq 4,\text{ }\left\Vert
\partial ^{\alpha }\varphi \right\Vert _{L^{\infty }\left( \left\vert
x\right\vert <2\right) }\leq C\text{ for }|\alpha |\leq 4, 
\]%
\[
\varphi (x)=|x|^{p}+O(|x|)\text{, for }|x|>1. 
\]

Let us consider the equation 
\begin{equation}
\partial _{t}\upsilon =i\left( \Delta \upsilon +A\upsilon +F\left(
x,t\right) \right) ,\text{ }x\in R^{n},\text{ }t\in \left[ 0,1\right] , 
\tag{4.10}
\end{equation}%
where $F\left( x,t\right) =\tilde{V}\upsilon ,$ $A$ is a symmetric operator
in $H$ and $\tilde{V}$ is a operator in $H$ defined by $\left( 3.3\right) .$

Let 
\[
f(x,t)=e^{\tilde{\gamma}\varphi }\upsilon (x,t),\text{ }Q\left( t\right)
=\left( f(x,t),f(x,t)\right) _{H}, 
\]%
where $\upsilon $ is a solution of $\left( 4.10\right) $. Then, by reasoning
as in Lemma 2.3 we have 
\begin{equation}
\partial _{t}f=Sf+Kf+i\left[ A+e^{\tilde{\gamma}\varphi }F\right] \text{, }%
\left( x,t\right) \in R^{n}\times \left[ 0,1\right] ,  \tag{4.11}
\end{equation}%
here $S$, $K$ are symmetric and skew-symmetric operator, respectively given
by%
\begin{equation}
S=-i\tilde{\gamma}\left( 2\nabla \varphi .\nabla +\Delta \varphi \right) 
\text{, }K=i\left( \Delta +A+\tilde{\gamma}^{2}\left\vert \nabla \varphi
\right\vert ^{2}\right) .  \tag{4.12}
\end{equation}%
Let 
\[
\left[ S,K\right] =SK-KS\text{.} 
\]%
A calculation shows that, 
\[
SK=\tilde{\gamma}\left( 2\nabla \varphi .\nabla +\Delta \varphi \right)
\left( \Delta +A+\tilde{\gamma}^{2}\left\vert \nabla \varphi \right\vert
^{2}\right) =\tilde{\gamma}\left( 2\nabla \varphi .\nabla +\Delta \varphi
\right) \Delta + 
\]%
\[
\tilde{\gamma}\left( 2\nabla \varphi .\nabla +\Delta \varphi \right) A+%
\tilde{\gamma}^{3}\left\vert \nabla \varphi \right\vert ^{2}\left( 2\nabla
\varphi .\nabla +\Delta \varphi \right) , 
\]%
\[
KS=\tilde{\gamma}\left[ \Delta \left( 2\nabla \varphi .\nabla +\Delta
\varphi \right) +A\left( 2\nabla \varphi .\nabla +\Delta \varphi \right) %
\right] +\tilde{\gamma}^{3}\left\vert \nabla \varphi \right\vert ^{2}\left(
2\nabla \varphi .\nabla +\Delta \varphi \right) , 
\]%
\[
\left[ S,K\right] =\tilde{\gamma}\left[ \left( 2\nabla \varphi .\nabla
+\Delta \varphi \right) \Delta -\Delta \left( 2\nabla \varphi .\nabla
+\Delta \varphi \right) \right] +2\tilde{\gamma}\left( \nabla \varphi
.\nabla A-A\nabla \varphi .\nabla \right) , 
\]%
\begin{equation}
\text{ }S_{t}+\left[ S,K\right] =-2i\tilde{\gamma}\left( \nabla \varphi
.\partial _{t}\nabla +\Delta \varphi \partial _{t}\right) +  \tag{4.13}
\end{equation}%
\[
\tilde{\gamma}\left[ \left( 2\nabla \varphi .\nabla +\Delta \varphi \right)
\Delta -\Delta \left( 2\nabla \varphi .\nabla +\Delta \varphi \right) \right]
+2\tilde{\gamma}\left( \nabla \varphi .\nabla A-A\nabla \varphi .\nabla
\right) . 
\]%
By Lemma 2.2 
\[
Q^{^{\prime \prime }}\left( t\right) =2\partial _{t}\func{Re}\left( \partial
_{t}f-Sf-Kf,f\right) _{X}+2\left( S_{t}f+\left[ S,K\right] f,f\right) _{X}+ 
\]

\begin{equation}
\left\Vert \partial _{t}f-Sf+Kf\right\Vert _{X}^{2}-\left\Vert \partial
_{t}f-Sf-Kf\right\Vert _{X}^{2},  \tag{4.14}
\end{equation}%
so, 
\begin{equation}
Q^{^{\prime \prime }}\left( t\right) \geq 2\partial _{t}\func{Re}\left(
\partial _{t}f-Sf-Kf,f\right) _{X}+2\left( S_{t}f+\left[ S,K\right]
f,f\right) _{X}.  \tag{4.15}
\end{equation}%
Multiplying $(4.15)$ by $t(1-t)$ and integrating in $t$ we obtain 
\[
\dint\limits_{0}^{1}t(1-t)\left( S_{t}f+\left[ S,K\right] f,f\right)
_{X}dt\leq 
\]%
\[
M_{0}\left[ \sup\limits_{t\in \left[ 0,1\right] }\left\Vert e^{\tilde{\gamma}%
\varphi }\upsilon \left( .,t\right) \right\Vert _{X}+\sup\limits_{t\in \left[
0,1\right] }\left\Vert e^{\tilde{\gamma}\varphi }F\left( .,t\right)
\right\Vert _{X}\right] . 
\]

This computation can be justified by parabolic regularization using the fact
that we already know the decay estimate for the solution of $\left(
4.10\right) $. Hence, combining $(3.8)$, $(4.3)$ and $(4.9$) it follows that 
\[
\tilde{\gamma}\dint\limits_{0}^{1}\dint\limits_{R^{n}}t(1-t)D^{2}\varphi
\left( x,t\right) \left( \nabla f,\nabla f\right) _{H}dxdt+\tilde{\gamma}%
^{3}\dint\limits_{0}^{1}\dint\limits_{R^{n}}t(1-t)D^{2}\varphi \left(
x,t\right) \left( \nabla f,\nabla f\right) _{H}dxdt\leq 
\]%
\begin{equation}
M_{0}\left[ \sup\limits_{t\in \left[ 0,1\right] }\left\Vert e^{\tilde{\gamma}%
\varphi }\upsilon \left( .,t\right) \right\Vert _{X}\left( 1+\left\Vert 
\tilde{V}\left( .,t\right) \right\Vert _{L^{\infty }\left( R^{n};L\left(
H\right) \right) }\right) +\tilde{\gamma}\sup\limits_{t\in \left[ 0,1\right]
}\left\Vert e^{\tilde{\gamma}\varphi }\upsilon \left( .,t\right) \right\Vert
_{X}\right] \leq  \tag{4.16}
\end{equation}%
\[
M_{0}k^{C_{p}}a_{k}. 
\]

It is clear to see that%
\[
\nabla f=\tilde{\gamma}\upsilon e^{\tilde{\gamma}\varphi }\nabla \varphi +e^{%
\tilde{\gamma}\varphi }\nabla \upsilon \text{.} 
\]%
So, by using the properties of $\varphi $ we get 
\[
\text{ }\left\vert e^{2\tilde{\gamma}\varphi }D^{2}\varphi \left\vert \nabla
\varphi \right\vert ^{2}\right\vert \leq C_{p}e^{3\gamma \varphi \mid 2}. 
\]

From here, we can conclude that 
\[
\gamma \dint\limits_{0}^{1}\dint\limits_{R^{n}}t(1-t)\left( 1+\left\vert
x\right\vert \right) ^{p-2}\left\Vert \nabla \upsilon \left( x,t\right)
\right\Vert _{H}e^{\gamma \left\vert x\right\vert ^{p}}dxdt+ 
\]%
\begin{equation}
\sup\limits_{t\in \left[ 0,1\right] }\left\Vert e^{\frac{\gamma \left\vert
x\right\vert ^{p}}{2}}\upsilon \left( .,t\right) \right\Vert _{X}\leq
C_{0}k^{C_{p}}a_{k}^{2}=C_{0}k^{C_{p}}e^{a\left( k,p\right) }  \tag{4.17}
\end{equation}%
for $k\geq k_{0}(M_{0})$ sufficiently large, where 
\[
a\left( k,p\right) =C_{\mu }M_{0}k^{C_{p}}e^{2a_{1}k^{1\mid \left(
2-p\right) }}. 
\]%
For proving Theorem 1 first, we deduce the following estimate%
\begin{equation}
\dint\limits_{\left\vert x\right\vert <\frac{R}{2}}\dint\limits_{\nu
_{1}}^{\nu _{2}}\left\Vert \tilde{u}\left( x,t\right) \right\Vert
_{H}dtdx\geq C_{\nu }e^{-M}\left\Vert u\left( .,0\right) \right\Vert _{X}, 
\tag{4.18}
\end{equation}%
for $r$ sufficiently large, $\nu _{1}$, $\nu _{2}\in \left( 0,1\right) $, $%
\nu _{1}<\nu _{2}$ and $\nu =\left( \nu _{1},\nu _{2}\right) $. From $\left(
3.1\right) $ by using the change of variables $s=\beta t\sigma \left(
t\right) $ and $y=\sqrt{\alpha \beta }x\sigma \left( t\right) $ we get 
\begin{equation}
\dint\limits_{\left\vert x\right\vert <\frac{r}{2}}\dint\limits_{\nu
_{1}}^{\nu _{2}}\left\Vert \tilde{u}\left( x,t\right) \right\Vert
_{H}^{2}dtdx=  \tag{4.19}
\end{equation}%
\[
\left( \alpha \beta \right) ^{\frac{n}{2}}\dint\limits_{\left\vert
x\right\vert <\frac{R}{2}}\dint\limits_{\nu _{1}}^{\nu _{2}}\left\vert
\sigma \left( t\right) \right\vert ^{n}\left\Vert u\left( \sqrt{\alpha \beta 
}x\sigma \left( t\right) ,\beta t\sigma \left( t\right) \right) \right\Vert
_{H}^{2}dtdx\geq 
\]%
\[
M_{0}\frac{\beta }{\alpha }\dint\limits_{\left\vert y\right\vert
<R_{0}}\dint\limits_{s\nu _{1}}^{s\nu _{2}}\left\Vert u\left( y,s\right)
\right\Vert _{H}^{2}\frac{dsdy}{s^{2}}\geq M_{0}\frac{\beta }{\alpha }%
\dint\limits_{\left\vert y\right\vert <R_{0}}\dint\limits_{s\nu _{1}}^{s\nu
_{2}}\left\Vert u\left( y,s\right) \right\Vert _{H}^{2}dsdy 
\]%
for $k>M_{0},$ $s\nu _{1}>\frac{1}{2}$ and $r_{0}=r\left( ka_{0}^{-1}\right)
^{\frac{1}{2p}}.$ Thus, taking 
\begin{equation}
r>\omega \left( ka_{0}^{-1}\right) ^{\frac{1}{2p}}  \tag{4.20}
\end{equation}%
with $\omega =\omega (u)$ a constant to be determined, it follows that%
\[
\Phi \geq M_{0}\frac{\beta }{\alpha }\dint\limits_{\left\vert y\right\vert
<\omega }\dint\limits_{s\nu _{1}}^{s\nu _{2}}\left\Vert u\left( y,s\right)
\right\Vert _{H}^{2}dsdy, 
\]%
where the interval $I=I_{k}=[s\nu _{1},s\nu _{2}]$ satisfies $I\subset
\lbrack 1/2,1]$ for $k$ sufficiently large. Moreover, given $\varepsilon >0$
there exists $k_{0}(\varepsilon )>0$ such that for any $k\geq k_{0}$ one has
that $I_{k}\subset \lbrack 1-\varepsilon ,1]$. By hypothesis on $u(x,t)$,
i.e. the continuity of $\left\Vert u(\text{\textperiodcentered }%
,s)\right\Vert _{X}$ at $s=1$, it follows that there exists $\omega >1$ and $%
K_{0}=K_{0}(u)$ such that for any $k\geq K_{0}$ and for any $s\in I_{k}$ 
\[
\dint\limits_{\left\vert y\right\vert <\omega }\left\Vert u\left( y,s\right)
\right\Vert _{H}^{2}dy\geq C_{\nu }e^{-M}\left\Vert u\left( .,0\right)
\right\Vert _{X}, 
\]%
which yields the desired result. Next, we deduce the following estimate%
\begin{equation}
\dint\limits_{\left\vert x\right\vert <R}\dint\limits_{\mu _{1}}^{\mu
_{2}}\left( \left\Vert \tilde{u}\left( x,t\right) \right\Vert
_{H}^{2}+\left\Vert \nabla \tilde{u}\left( x,t\right) \right\Vert
_{H}^{2}+\left\Vert A\tilde{u}\left( x,t\right) \right\Vert _{H}^{2}\right)
dtdx\leq  \tag{4.21}
\end{equation}%
\[
C_{n\mu }C_{0}k^{C_{p}}e^{a\left( k,p\right) }, 
\]%
for $r$ sufficiently large, $\nu _{1},$ $\nu _{2}\in \left( 0,1\right) $, $%
\mu _{1}=\frac{\left( \nu _{2}-\nu _{1}\right) }{8}$, $\mu _{2}=1-\mu _{1},$ 
$\mu _{1}<\mu _{2}$ and $\mu =\left( \mu _{1},\mu _{2}\right) $. Indeed,
from $\left( 3.9\right) $ and $\left( 4.17\right) $\ we obtain 
\[
\dint\limits_{\left\vert x\right\vert <R}\dint\limits_{\mu _{1}}^{\mu
_{2}}\left\Vert \tilde{u}\left( x,t\right) \right\Vert _{H}dtdx\leq C_{\mu
}e^{2M}\left\Vert u\left( .,0\right) \right\Vert _{X}, 
\]%
\begin{equation}
\dint\limits_{\mu _{1}}^{\mu _{2}}\dint\limits_{\left\vert x\right\vert
<R}\left\Vert \nabla \tilde{u}\left( x,t\right) \right\Vert _{H}dtdx\leq
C_{n\mu }\dint\limits_{\mu _{1}}^{\mu _{2}}\dint\limits_{\left\vert
x\right\vert <R}t(1-t)\left\Vert \nabla \upsilon \left( x,t\right)
\right\Vert _{H}e^{\gamma \left\vert x\right\vert ^{p}}dtdx\leq  \tag{4.22}
\end{equation}%
\[
C_{\mu }\gamma ^{-1}r^{2-p}C_{0}k^{C_{p}}A_{k}^{2}\leq C_{\mu
}C_{0}k^{C_{p}}e^{2a_{1}k^{\left( 2-p\right) ^{-1}}}. 
\]

Hence, from $\left( 4.22\right) $ we get $\left( 4.21\right) $ for $k\geq
k_{0}(C_{0})$ sufficiently large.

Let $Y=L^{2}\left( R^{n}\times \left[ 0,1\right] ;H\right) $. By reasoning
as in $\left[ 6\text{, Lemma 3.1}\right] $ we obtain

\textbf{Lemma 4.1}. Assume the assumpt\i ions (1) and \ (3) of Condition 1
are satisfied. Suppose that $r>0$ and $\varphi $ : $[0,1]\rightarrow \mathbb{%
R}$ is a smooth function. Then, there exists $C=C(n,\varphi ,H,A)>0$ such
that, the inequality%
\[
\dsum\limits_{k=1}^{n}x_{k}\left[ \left( \frac{\partial A}{\partial x_{k}}%
g,g\right) +\left( \frac{\partial g}{\partial x_{k}},Ag\right) _{X}\right] +%
\frac{\varkappa ^{\frac{3}{2}}}{r^{2}}\left\Vert e^{\varkappa \left\vert
\psi \right\vert }g\right\Vert _{Y}\leq C\left\Vert e^{\varkappa \left\vert
\psi \right\vert }i\left( \partial _{t}g+\Delta g+Ag\right) \right\Vert _{Y} 
\]%
holds, for $\varkappa \geq CR^{2}$ and $g\in C_{0}^{\infty }\left(
R^{n+1};H\left( A\right) \right) $ with support contained in the set 
\[
\left\{ x,t:\text{ }\left\vert \psi \left( x,t\right) \right\vert
=\left\vert \frac{x}{r}+\varphi \left( t\right) e_{1}\right\vert \geq
1\right\} . 
\]

\textbf{Proof. }Let $f=e^{\alpha \left\vert \psi \left( x,t\right)
\right\vert ^{2}}g.$ Then, by acts of Schredinger operator $\left( i\partial
_{t}+\Delta +A\right) $ to $f\in X\left( A\right) $ we get 
\[
e^{\alpha \left\vert \psi \left( x,t\right) \right\vert ^{2}}\left( \left(
i\partial _{t}g+\Delta g+Ag\right) \right) =S_{\alpha }f-4\alpha K_{\alpha
}f, 
\]%
where 
\[
S_{\alpha }=\left( i\partial _{t}+\Delta +A\right) +\frac{4\alpha ^{2}}{r^{2}%
}\text{ }\left\vert \psi \left( x,t\right) \right\vert ^{2},\text{ } 
\]%
\[
K_{\alpha }=\frac{1}{r}\left( \frac{x}{r}+\varphi \left( t\right)
e_{1}\right) .\nabla +\frac{n}{r^{2}}+\frac{i\varphi ^{\prime }}{2}\left( 
\frac{x_{1}}{r}+\varphi \left( t\right) \right) . 
\]%
Hence, 
\[
\left( S_{\alpha }\right) ^{\ast }=S_{\alpha }\text{, }\left( K_{\alpha
}\right) ^{\ast }=K_{\alpha } 
\]%
and 
\[
\left\Vert e^{\alpha \left\vert \psi \left( x,t\right) \right\vert
^{2}}\left( i\partial _{t}g+\Delta g+Ag\right) \right\Vert _{X}^{2}=\left(
S_{\alpha }f-4\alpha K_{\alpha }f,S_{\alpha }f-4\alpha K_{\alpha }f\right)
_{X}\geq 
\]%
\[
-4\alpha \left( \left( S_{\alpha }K_{\alpha }-K_{\alpha }S_{\alpha }\right)
f,f\right) _{X}=-4\alpha \left( \left[ S_{\alpha },K_{\alpha }\right]
f,f\right) _{X}. 
\]%
A calculation shows that 
\begin{equation}
\left[ S_{\alpha },K_{\alpha }\right] =\frac{2}{r^{2}}\Delta -\frac{4\alpha
^{2}}{r^{4}}\left\vert \frac{x}{r}+\varphi e_{1}\right\vert ^{2}-\frac{1}{2}+%
\frac{2i\varphi ^{\prime }}{r}\partial _{x_{1}}+\left[ A,K_{\alpha }\right] ,
\tag{4.23}
\end{equation}%
where%
\[
\left[ A,K_{\alpha }\right] f=\left[ \left( \frac{x_{1}}{r}+\varphi \right) A%
\frac{\partial f}{\partial x_{1}}+\frac{1}{r}\dsum\limits_{k=2}^{n}x_{k}A%
\frac{\partial f}{\partial x_{k}}\right] - 
\]%
\[
\left( \frac{x_{1}}{r}+\varphi \right) \left( \frac{\partial A}{\partial
x_{1}}f+A\frac{\partial f}{\partial x_{1}}\right) +\frac{1}{r}%
\dsum\limits_{k=2}^{n}x_{k}\left( \frac{\partial A}{\partial x_{k}}f+A\frac{%
\partial f}{\partial x_{k}}\right) . 
\]%
Since $A$ is a symmetric operator in $H,$\ from $\left( 4.23\right) $ we
have 
\begin{equation}
\left\Vert e^{\alpha \left\vert \psi \left( x,t\right) \right\vert
^{2}}\left( i\partial _{t}g+\Delta g+Ag\right) \right\Vert _{X}^{2}\geq 
\tag{4.24}
\end{equation}%
\[
\frac{16\alpha ^{3}}{r^{4}}\dint \left\vert \frac{x}{r}+\varphi
e_{1}\right\vert ^{2}\left\Vert f\left( x,t\right) \right\Vert ^{2}dxdt+%
\frac{8\alpha }{r^{2}}\dint \left\Vert \nabla f\left( x,t\right) \right\Vert
^{2}dxdt+ 
\]%
\[
2\alpha \dint \left[ \left( \frac{x_{1}}{r}+\varphi \right) \varphi ^{\prime
\prime }+\left( \varphi ^{\prime }\right) ^{2}\right] \left\Vert f\left(
x,t\right) \right\Vert ^{2}dxdt-\frac{8\alpha }{r}\dint \varphi ^{\prime
}\left( \partial _{x_{11}}f,\bar{f}\right) dxdt+ 
\]%
\[
C\dsum\limits_{k=1}^{n}x_{k}\left[ \left( \frac{\partial A}{\partial x_{k}}%
f,f\right) +\left( \frac{\partial f}{\partial x_{k}},Af\right) _{X}\right] . 
\]

In view of the hypothesis on $g$ and the Cauchy--Schwarz inequality, the
absolute value of the third fourth terms in $(4.23)$ can be bounded by a
fraction of the first two terms on the right-hand side of $(4.24)$, when $%
\alpha >Cr^{2}$ for some large $C$ depending on $\left\Vert \varphi ^{\prime
}\right\Vert _{\infty }+\left\Vert \varphi ^{\prime \prime }\right\Vert
_{\infty }.$ Moreover, by using\ the assumption on $A=A\left( x\right) $, we
get that the last two terms are nonnegative. This yields the assertion. Now,
from $\left( 3.3\right) $ we have \ \ \ \ \ 
\[
\left\Vert \tilde{V}\left( x,t\right) \right\Vert _{H}\leq \frac{\alpha }{%
\beta }\mu _{1}^{-2}\left\Vert V\right\Vert _{B}\leq \mu _{1}^{-2}a_{0}^{%
\frac{1}{p}}k^{-\frac{1}{p}}\left\Vert V\right\Vert _{B}. 
\]%
Then from $\left( 4.20\right) $ we get 
\begin{equation}
\left\Vert \tilde{V}\right\Vert _{L^{\infty }\left( R^{n}\times \left[ \mu
_{1},\mu _{2}\right] ;L\left( H\right) \right) }<R.  \tag{4.25}
\end{equation}%
Define 
\[
\delta ^{2}\left( r\right) = 
\]

\[
\dint\limits_{\mu _{1}}^{\mu _{2}}\dint\limits_{R-1\leq \left\vert
x\right\vert \leq r}\left( \left\Vert \tilde{u}\left( x,t\right) \right\Vert
_{H}^{2}+\left\Vert \nabla \tilde{u}\left( x,t\right) \right\Vert
_{H}^{2}+\left\Vert A\tilde{u}\left( x,t\right) \right\Vert _{H}^{2}\right)
dtdx. 
\]%
Let $\nu _{1}$, $\nu _{2}\in \left( 0,1\right) $, $\nu _{1}<\nu _{2}$ and $%
\nu _{2}<2\nu _{1}$. We choose $\varphi \in C^{\infty }\left( \left[ 0,1%
\right] \right) $ and $\theta $, $\theta _{r}\in C_{0}^{\infty }\left(
R^{n}\right) $ satisfying%
\[
0\leq \varphi \left( t\right) \leq 3\text{, }\varphi \left( t\right) =3\text{
for }t\in \left[ \nu _{1},\nu _{2}\right] \text{, }\varphi \left( t\right) =0%
\text{ for } 
\]%
\[
t\in \left[ 0,\nu _{2}-\nu _{1}\right] \cup \left[ \nu _{2}+\frac{\nu
_{2}-\nu _{1}}{2},1\right] , 
\]

\[
\theta _{r}\left( x\right) =1\text{ for }\left\vert x\right\vert <r-1,\text{ 
}\theta _{R}\left( x\right) =0\text{, for }\left\vert x\right\vert >r\text{,}
\]%
and 
\[
\theta \left( x\right) =1\text{ for }\left\vert x\right\vert <1,\text{ }%
\theta \left( x\right) =0\text{, for }\left\vert x\right\vert \geq 2\text{.} 
\]

Let 
\begin{equation}
g\left( x,t\right) =\theta _{r}\left( x\right) \theta \left( \psi \left(
x,t\right) \right) \tilde{u}\left( x,t\right) ,  \tag{4.26}
\end{equation}%
where $\tilde{u}\left( x,t\right) $ is a solution of $\left( 3.2\right) $
when $\tilde{V}=\tilde{F}\equiv 0$. It is clear to see that 
\[
\left\vert \psi \left( x,t\right) \right\vert \geq \frac{5}{2}\text{ for }%
\left\vert x\right\vert <\frac{r}{2}\text{ and }t\in \left[ \nu _{1},\nu _{2}%
\right] . 
\]%
Hence, 
\[
g\left( x,t\right) =\tilde{u}\left( x,t\right) \text{ and }e^{\varkappa
\left\vert \psi \left( x,t\right) \right\vert ^{2}}\geq e^{\frac{25}{4}%
\varkappa }\text{ for }\left\vert x\right\vert <\frac{r}{2}\text{, }t\in %
\left[ \nu _{1},\nu _{2}\right] . 
\]%
Moreover, from $\left( 4.26\right) $ also we get that 
\[
g\left( x,t\right) =0\text{ for }\left\vert x\right\vert \geq r\text{ or }%
t\in \left[ 0,\nu _{2}-\nu _{1}\right] \cup \left[ \nu _{2}+\frac{\nu
_{2}-\nu _{1}}{2},1\right] , 
\]%
so 
\[
\text{supp }g\subset \left\{ \left\vert x\right\vert \leq R\right\} \times %
\left[ \nu _{2}-\nu _{1},\nu _{2}+\frac{\nu _{2}-\nu _{1}}{2}\right] \cap
\left\{ \left\vert b\left( x,t\right) \right\vert \geq 1\right\} . 
\]%
Then, for $\xi =\psi \left( x,t\right) $ we have%
\[
\left( i\partial _{t}+\Delta +A+\tilde{V}\right) g=\left[ \theta \left( \xi
\right) \left( 2\nabla \theta _{r}\left( x\right) .\tilde{u}+\tilde{u}\Delta
\theta _{r}\left( x\right) \right) +2\nabla \theta \left( \xi \right)
.\nabla \theta _{r}\tilde{u}\right] + 
\]%
\[
\theta _{r}\left( x\right) \left[ 2r^{-1}\nabla \theta \left( \xi \right)
.\nabla \tilde{u}+r^{-2}\tilde{u}\Delta \theta \left( \xi \right) +i\varphi
^{\prime }\partial x_{1}\theta \left( \xi \right) u\right] =B_{1}+B_{2}\text{%
.} 
\]%
Note that, 
\[
\text{supp }B_{1}\subset \left\{ \left( x,t\right) :\text{ }r-1\leq
\left\vert x\right\vert \leq r,\text{ }\mu _{1}\leq t\leq \mu _{2}\right\} 
\]%
and 
\[
\text{supp }B_{2}\subset \left\{ \left( x,t\right) \in R^{n}\times \left[ 0,1%
\right] \text{, }1\leq \left\vert \psi \left( x,t\right) \right\vert \leq
2\right\} . 
\]%
Now applying Lemma 4.1 choosing $\varkappa =d_{n}^{2}R^{2},$ $d_{n}^{2}\geq
\left\Vert \varphi ^{\prime }\right\Vert _{\infty }+\left\Vert \varphi
^{\prime \prime }\right\Vert _{\infty }$ it follows that

\begin{equation}
R\left\Vert e^{\varkappa \left\vert \psi \right\vert ^{2}}g\right\Vert
_{Y}\leq C\left\Vert e^{\varkappa \left\vert \psi \right\vert ^{2}}i\left(
\partial _{t}g+\Delta g+Ag\right) \right\Vert _{Y}\leq  \tag{4.27}
\end{equation}

\[
C\left[ \left\Vert e^{\varkappa \left\vert \psi \right\vert ^{2}}\tilde{V}%
g\right\Vert _{Y}+\left\Vert e^{\varkappa \left\vert \psi \right\vert
^{2}}B_{1}\right\Vert _{Y}+\left\Vert e^{\varkappa \left\vert \psi
\right\vert ^{2}}B_{2}\right\Vert _{Y}\right] = 
\]%
\[
D_{1}+D_{2}+D_{3}. 
\]%
Since 
\[
\left\Vert \tilde{V}\right\Vert _{L^{\infty }\left( R^{n}\times \left[ \mu
_{1},\mu _{2}\right] ;L\left( H\right) \right) }<r, 
\]%
$D_{1}$ can be absorbed in the left hand side of $\left( 4.27\right) .$
Moreover, $\left\vert \psi \left( x,t\right) \right\vert \leq 4$ on the
support of $B_{1},$ thus \ 
\[
D_{2}\leq C\delta \left( R\right) e^{16\varkappa }. 
\]%
Let%
\[
R_{\mu }^{n}=\left\{ \left( x,t\right) :\left\vert x\right\vert \leq r,\text{
}\mu _{1}\leq t\leq \mu _{2}\right\} 
\]%
Then $R_{\mu }^{n}\subset $supp $B_{2}$, and $1\leq \left\vert \psi \left(
x,t\right) \right\vert \leq 2$, so%
\[
D_{3}\leq Ce^{4\varkappa }\left\Vert \tilde{u}+\nabla \tilde{u}\right\Vert
_{L^{2}\left( R^{n}\times \left[ \mu _{1},\mu _{2}\right] ;H\right) }. 
\]%
By using $\left( 4.18\right) $ and $\left( 4.22\right) $ we have 
\begin{equation}
C_{\mu }e^{-M}e^{\frac{25}{4}\varkappa }\left\Vert u\left( .,0\right)
\right\Vert _{X}\leq re^{\frac{25}{4}\varkappa }\left[ \dint\limits_{\left%
\vert x\right\vert <\frac{R}{2}}\dint\limits_{\nu _{1}}^{\nu _{2}}\left\Vert 
\tilde{u}\left( x,t\right) \right\Vert _{H}dtdx\right] ^{\frac{1}{2}}\leq 
\tag{4.28}
\end{equation}%
\[
C_{\mu }\delta \left( r\right) e^{\frac{25}{4}\varkappa }+C_{\mu
}e^{4\varkappa }\left\Vert \tilde{u}+\nabla \tilde{u}\right\Vert
_{L^{2}\left( R_{\mu }^{n};H\right) }\leq C_{\mu }\delta \left( r\right)
e^{16\varkappa }+C_{\mu }C_{0}k^{C_{p}}e^{4\varkappa }e^{2a_{1}k^{\left(
2-p\right) ^{-1}}}. 
\]%
Puting $\varkappa =d_{n}r^{2}=2a_{1}k^{\frac{1}{2-p}}$ it follows from $%
\left( 4.28\right) $ that, if $\left\Vert u\left( .,0\right) \right\Vert
_{X}\neq 0$ then%
\begin{equation}
\delta \left( r\right) \geq C_{\mu }\left\Vert u\left( .,0\right)
\right\Vert _{X}e^{-\left( M+10\varkappa \right) }=C_{\mu }\left\Vert
u\left( .,0\right) \right\Vert _{X}e^{-\left( M+20\right) a_{1}k^{\frac{1}{%
2-p}}}  \tag{4.29 }
\end{equation}%
for $k\geq k_{0}(C_{\mu })$ sufficiently large. Now, by $\left( 4.22\right) $
we get%
\[
\delta ^{2}\left( r\right) =\dint\limits_{\mu _{1}}^{\mu
_{2}}\dint\limits_{r-1\leq \left\vert x\right\vert \leq r}\left( \left\Vert 
\tilde{u}\left( x,t\right) \right\Vert _{H}^{2}+\left\Vert \nabla \tilde{u}%
\left( x,t\right) \right\Vert _{H}^{2}+\left\Vert A\tilde{u}\left(
x,t\right) \right\Vert _{H}^{2}\right) dtdx\leq 
\]%
\[
\dint\limits_{\mu _{1}}^{\mu _{2}}\dint\limits_{r-1\leq \left\vert
x\right\vert \leq r}\left( \left\Vert \tilde{u}\left( x,t\right) \right\Vert
_{H}^{2}+\left\Vert A\tilde{u}\left( x,t\right) \right\Vert _{H}^{2}\right)
dtdx+ 
\]%
\[
C_{\mu }\dint\limits_{\mu _{1}}^{\mu _{2}}\dint\limits_{r-1\leq \left\vert
x\right\vert \leq r}t\left( 1-t\right) \left\Vert \nabla \tilde{u}\left(
x,t\right) \right\Vert _{H}^{2}dtdx\leq 
\]%
\[
C_{\mu }e^{-\gamma \left( r-1\right) ^{p}}\sup\limits_{t\in \left[ 0,1\right]
}\left\Vert e^{\gamma \left\vert x\right\vert ^{p/2}}\tilde{u}\left(
x,t\right) \right\Vert _{X}^{2}+C_{\mu }\gamma ^{-1}r^{2-p}e^{-\gamma \left(
r-1\right) ^{p}}\times 
\]%
\begin{equation}
\dint\limits_{\mu _{1}}^{\mu _{2}}\dint\limits_{r-1\leq \left\vert
x\right\vert \leq r}\left[ \frac{t\left( 1-t\right) }{\left( 1+\left\vert
x\right\vert \right) ^{2/p}}\left\Vert \nabla \tilde{u}\left( x,t\right)
\right\Vert _{H}^{2}+\left\Vert A\tilde{u}\left( x,t\right) \right\Vert
_{H}^{2}\right] dtdx\leq  \tag{4.30}
\end{equation}%
\[
C_{\mu }\gamma ^{-1}k^{C_{p}}e^{\eta \left( p\right) },\text{ }\eta \left(
p\right) =2a_{1}k^{\left( 2-p\right) ^{-1}}-\gamma \left( r-1\right) ^{p}. 
\]%
The estimates $\left( 4.28\right) -\left( 4.30\right) $ imply 
\begin{equation}
C_{\mu }e^{-2M}e^{\frac{25}{4}\varkappa }\left\Vert u\left( .,0\right)
\right\Vert _{X}\leq C_{0}k^{C_{p}}e^{\omega \left( p\right) }+O\left(
k^{1/2\left( 2-p\right) }\right) \text{,}  \tag{4.31}
\end{equation}%
where 
\[
\omega \left( p\right) =42a_{1}k^{1/\left( 2-p\right) }-a_{0}^{-\frac{1}{2}%
}\left( \frac{2a_{1}}{d_{n}}\right) ^{\frac{p}{2}}k^{1/\left( 2-p\right) }. 
\]%
Hence, if $42a_{1}<\sqrt{a_{1}^{p}a_{0}}\left( \frac{2}{d_{n}}\right) ^{%
\frac{p}{2}}$ by letting $k$ tends to infinity it follows from $\left(
4.31\right) $ that $\left\Vert u\left( .,0\right) \right\Vert _{X}=0$, which
gives $u\left( x,t\right) \equiv 0$.

\textbf{Proof of Corollary 1. }Since 
\[
\dint\limits_{R^{n}}\left\Vert u\left( x,1\right) \right\Vert
_{H}^{2}e^{2b\left\vert x\right\vert ^{q}}dx<\infty \text{ for }b=\frac{%
\beta ^{q}}{q} 
\]%
one has that%
\[
\dint\limits_{R^{n}}\left\Vert u\left( x,1\right) \right\Vert
_{H}^{2}e^{2k\left\vert x\right\vert ^{q}}dx\leq \left\Vert e^{2k\left\vert
x\right\vert ^{q}-2b\left\vert x\right\vert ^{q}}\right\Vert _{\infty
}\dint\limits_{R^{n}}\left\Vert u\left( x,1\right) \right\Vert
_{H}^{2}e^{2b\left\vert x\right\vert ^{q}}dx. 
\]

Then, by reasoning as in $\left[ \text{7, Corollary 1}\right] $ we obtain
the assertion.

\textbf{\ Proof of Theorem 2. }Indeed, just applying Corollary 1 with 
\[
u\left( x,t\right) =u_{1}\left( x,t\right) -u_{2}\left( x,t\right) 
\]%
and 
\[
V\left( x,t\right) =\frac{F\left( u_{1},\bar{u}_{1}\right) -F\left( u_{2},%
\bar{u}_{2}\right) }{u_{1}-u_{2}} 
\]%
we obtain the assertion of Theorem 2.

\begin{center}
\textbf{5. Proof of Theorem 3. }
\end{center}

First, we deduce the corresponding upper bounds. Assume 
\[
\left\Vert u\left( .,t\right) \right\Vert _{X}=a\neq 0. 
\]

Fix $\bar{t}\in \left( 0,1\right) $ near $1,$ and let 
\[
\upsilon \left( x,t\right) =u\left( x,t-1+\bar{t}\right) \text{, }t\in \left[
0,1\right] 
\]%
which satisfies the equation $\left( 2.10\right) $ with%
\begin{equation}
\left\vert \upsilon \left( x,0\right) \right\vert \leq \frac{b_{1}}{\left( 2-%
\bar{t}\right) ^{n/2}}e^{-\frac{b_{2}\left\vert x\right\vert ^{p}}{\left( 2-%
\bar{t}\right) ^{p}}},\text{ }\left\vert \upsilon \left( x,1\right)
\right\vert \leq \frac{b_{1}}{\left( 1-\bar{t}\right) ^{n/2}}e^{-\frac{%
b_{2}\left\vert x\right\vert ^{p}}{\left( 1-\bar{t}\right) ^{p}}}  \tag{5.1}
\end{equation}%
where $A$ is a linear operator$,$ $V\left( x,t\right) $ is a given potential
operator function in a Hilbert space $H.$

From $\left( 5.1\right) $ we get%
\[
\dint\limits_{R^{n}}\left\Vert \upsilon \left( x,0\right) \right\Vert
_{H}^{2}e^{A_{0}\left\vert x\right\vert ^{q}}dx=a_{0}^{2},\text{ }%
\dint\limits_{R^{n}}\left\Vert \upsilon \left( x,1\right) \right\Vert
_{H}^{2}e^{A_{1}\left\vert x\right\vert ^{q}}dx=A_{1}^{2}, 
\]%
where 
\begin{equation}
A_{0}=\frac{b_{2}}{\left( 2-\bar{t}\right) ^{p}}\text{, }A_{1}=\frac{b_{2}}{%
\left( 1-\bar{t}\right) ^{p}}.  \tag{5.2}
\end{equation}

For $V\left( x,t\right) =F\left( u,\bar{u}\right) $, by hypothesis 
\[
\left\Vert V\left( x,t\right) \right\Vert _{H}\leq C\left\Vert u\left( x,t-1+%
\bar{t}\right) \right\Vert _{H}^{\theta }\leq \frac{C}{\left( 2-t-\bar{t}%
\right) ^{\theta n/2}}e^{\frac{C\left\vert x\right\vert ^{p}}{\left( 2-t-%
\bar{t}\right) ^{p}}}. 
\]%
By using Appell transformation if we suppose that $\upsilon \left(
y,s\right) $ is a solution of 
\[
\partial _{s}\upsilon =i\left( \Delta \upsilon +A\upsilon +V\left(
y,s\right) \upsilon \right) \text{, }y\in R^{n},\text{ }t\in \left[ 0,1%
\right] , 
\]%
$\alpha $ and $\beta $ are positive, then 
\begin{equation}
\tilde{u}\left( x,t\right) =\left( \sqrt{\alpha \beta }\sigma \left(
t\right) \right) ^{\frac{n}{2}}u\left( \sqrt{\alpha \beta }x\sigma \left(
t\right) ,\beta t\sigma \left( t\right) \right) e^{\eta }.  \tag{5.3}
\end{equation}%
\ verifies the equation 
\begin{equation}
\partial _{t}\tilde{u}=i\left[ \Delta \tilde{u}+A\tilde{u}+\tilde{V}\left(
x,t\right) \tilde{u}+\tilde{F}\left( x,t\right) \right] ,\text{ }x\in R^{n},%
\text{ }t\in \left[ 0,1\right]  \tag{5.4}
\end{equation}%
with $\tilde{V}\left( x,t\right) ,$ $\tilde{F}\left( x,t\right) $ defined by 
$\left( 3.3\right) $, $\left( 3.4\right) $ and 
\begin{equation}
\left\Vert e^{\gamma \left\vert x\right\vert ^{p}}\tilde{u}_{k}\left(
x,0\right) \right\Vert _{X}=\left\Vert e^{\gamma \left( \frac{\alpha }{\beta 
}\right) ^{p/2}\left\vert x\right\vert ^{p}}\upsilon \left( x,0\right)
\right\Vert _{X}=a_{0},  \tag{5.5}
\end{equation}

\[
\left\Vert e^{\gamma \left\vert x\right\vert ^{p}}\tilde{u}_{k}\left(
x,1\right) \right\Vert _{X}=\left\Vert e^{\gamma \left( \frac{\beta }{\alpha 
}\right) ^{p/2}\left\vert x\right\vert ^{p}}\upsilon \left( x,1\right)
\right\Vert _{X}=a_{1}. 
\]%
\ It follows from expressions $\left( 3.6\right) $ and $\left( 3.8\right) $\
that 
\begin{equation}
\gamma \sim \frac{1}{\left( 1-\bar{t}\right) ^{p/2}}\text{, }\beta \sim 
\frac{1}{\left( 1-\bar{t}\right) ^{p/2}}\text{, }\alpha \sim 1.  \tag{5.6}
\end{equation}%
Next, we shall estimate 
\[
\left\Vert \tilde{V}\left( x,t\right) \right\Vert _{L_{t}^{1}L_{x}^{\infty
}\left( L\left( H\right) ,R\right) },\text{ } 
\]%
for a $r>0$, where 
\[
\text{ }L_{t}^{1}L_{x}^{\infty }\left( L\left( H\right) ,r\right)
=L^{1}\left( 0,1;L^{\infty }\left( R^{n}/O_{r}\right) ;L\left( H\right)
\right) . 
\]%
Thus, 
\[
\left\Vert \tilde{V}\left( x,t\right) \right\Vert _{L\left( H\right) }\leq 
\frac{\beta }{\alpha }\left\Vert V\left( y,s\right) \right\Vert _{L\left(
H\right) }\leq \frac{\beta }{\alpha }\frac{C}{\left( 1-\bar{t}\right)
^{\theta n/2}}e^{C\left\vert y\right\vert ^{p}}, 
\]%
with 
\[
\left\vert y\right\vert =\sqrt{\alpha \beta }\left\vert x\right\vert \sigma
\left( t\right) \geq R\sqrt{\frac{\alpha }{\beta }\sim }\frac{R}{\sqrt{\beta 
}}=Cr\left( 1-\bar{t}\right) ^{1/2}. 
\]%
Hence, 
\[
\left\Vert \tilde{V}\left( .,t\right) \right\Vert _{L^{\infty }\left(
R^{n};L\left( H\right) \right) }\leq \frac{\beta }{\alpha }\left\Vert
V\left( .,s\right) \right\Vert _{L^{\infty }\left( R^{n};L\left( H\right)
\right) }\leq \frac{C}{\left( 1-\bar{t}\right) ^{1+\theta n/2}} 
\]%
and 
\begin{equation}
\left\Vert \tilde{V}\left( x,t\right) \right\Vert _{L_{t}^{1}L_{x}^{\infty
}\left( L\left( H\right) ,r\right) }\leq \frac{\beta }{\alpha }\left\Vert
V\left( y,s\right) \right\Vert _{L_{t}^{1}L_{x}^{\infty }\left( r;\beta
\right) }\leq  \tag{5.7}
\end{equation}%
\[
\frac{C}{\left( 1-\bar{t}\right) ^{1+\theta n/2}}e^{-Cr^{p}\left( 1-\bar{t}%
\right) ^{1/2}}, 
\]%
where 
\[
L_{t}^{1}L_{x}^{\infty }\left( r,\beta \right) =L^{1}\left( 0,1;L^{\infty
}\left( R^{n}/O_{Cr/\sqrt{\beta }};L\left( H\right) \right) \right) . 
\]%
To apply Lemma 3.1 we need 
\begin{equation}
\left\Vert \tilde{V}\left( x,t\right) \right\Vert _{L_{t}^{1}L_{x}^{\infty
}\left( L\left( H\right) ,R\right) }\leq \frac{C}{\left( 1-\bar{t}\right)
^{1+\theta n/2}}e^{-Cr^{p}\left( 1-\bar{t}\right) ^{1/2}}\leq \delta _{0}. 
\tag{5.8}
\end{equation}%
for some $R,$ i.e., 
\begin{equation}
r\sim \frac{C}{\left( 1-\bar{t}\right) ^{p/2}}\delta \left( t\right) \text{, 
}  \tag{5.9}
\end{equation}%
where 
\[
\delta \left( t\right) =\log ^{\frac{1}{p}}\phi \left( t\right) \text{, }%
\phi \left( t\right) =\frac{C}{\delta _{0}\left( 1-\bar{t}\right) ^{\theta
n/2}}. 
\]%
Let 
\begin{equation}
\mathbb{V=}\tilde{V}_{\chi \left( \left\vert x>R\right\vert \right) }\left(
x,t\right) ,\text{ }\mathbb{F}=\tilde{V}_{\chi \left( \left\vert
x<r\right\vert \right) }\left( x,t\right) \tilde{u}\left( x,t\right) . 
\tag{5.10}
\end{equation}%
By using $\left( 5.5\right) -\left( 5.10\right) ,$ by virtue of Lemma 3.1
and $\left( 4.7\right) $ we deduced 
\[
\sup\limits_{t\in \left[ 0,1\right] }\left\Vert e^{\gamma \left\vert
x\right\vert ^{p}}\tilde{u}\left( .,t\right) \right\Vert _{X}^{2}\leq
C\left( \left\Vert e^{\gamma \left\vert x\right\vert ^{p}}\tilde{u}\left(
.,0\right) \right\Vert _{X}^{2}+\left\Vert e^{\gamma \left\vert x\right\vert
^{p}}\tilde{u}\left( .,1\right) \right\Vert _{X}^{2}\right) + 
\]%
\[
Ca^{2}e^{C\gamma r^{p}}\left\Vert \tilde{V}\left( x,t\right) \right\Vert
\leq \frac{Ca^{2}}{\left( 1-\bar{t}\right) ^{p}}e^{\delta \left( t\right) },%
\text{ }a=\left\Vert u_{0}\right\Vert _{X}. 
\]%
Next, using the same argument given in section $4$, $(4.10)-(4.22)$, one
finds that 
\[
\gamma \dint\limits_{0}^{1}\dint\limits_{R^{n}}t(1-t)\left( 1+\left\vert
x\right\vert \right) ^{p-2}\left\Vert \nabla \tilde{u}\left( x,t\right)
\right\Vert _{H}e^{\gamma \left\vert x\right\vert ^{p}}dxdt\leq \frac{Ca^{2}%
}{\left( 1-\bar{t}\right) ^{p}}e^{\delta \left( t\right) }. 
\]

Now we turn to the lower bounds estimates. Since they are similar to those
given in detail in Section 3, we obtain that the estimate $(4.24)$ for
potential operator function $\tilde{V}\left( x,t\right) $ when 
\[
\frac{\theta n}{2}-1<\frac{p}{2\left( 2-p\right) }\text{, i.e., }p>\frac{%
2\left( \theta n-2\right) }{\theta n-1}. 
\]%
Finally, we get 
\[
e^{\frac{C}{\left( 1-\bar{t}\right) ^{p}}\log \phi \left( t\right) }\leq
e^{C\gamma r^{p}}\leq e^{\frac{C}{\left( 1-\bar{t}\right) ^{p/2}}\vartheta
\left( t\right) },\text{ }\vartheta \left( t\right) =C\left( 1-\bar{t}%
\right) ^{\frac{-p^{2}}{2\left( 2-p\right) }} 
\]%
for $p<\frac{p}{2}+\frac{p^{2}}{2\left( 2-p\right) }$, i.e. $p>1$ that
assumed, i.e. we obtain the assertion of Theorem 3.

\textbf{Remark 5.1. }Let us consider the case $\theta =4/n$ in Theorem $3$,
i.e. $p>4/3$. Then from Theorem 3 we obtain

\textbf{Result 5.1}. Assume that the conditions of Theorem 3 are satisfied
for $p>4/3$. Then $u\left( x,t\right) \equiv 0.$

\begin{center}
\textbf{6. Unique continuation properties for the system of Schredinger
equation }
\end{center}

Consider the Cauchy problem for the finite or infinite system of Schr\"{o}%
dinger equation%
\begin{equation}
\frac{\partial u_{m}}{\partial t}=i\left[ \Delta
u_{m}+\sum\limits_{j=1}^{N}a_{mj}\left( x\right)
u_{j}+\sum\limits_{j=1}^{N}b_{mj}\left( x\right) u_{j}\right] ,\text{ }x\in
R^{n},\text{ }t\in \left( 0,T\right) ,  \tag{6.1}
\end{equation}%
where $u=\left( u_{1},u_{2},...,u_{N}\right) ,$ $u_{j}=u_{j}\left(
x,t\right) ,$ $a_{mj}=a_{mj}\left( x\right) $ are and $b_{mj}=b_{mj}\left(
x,t\right) $ are complex valued functions$.$ Let $l_{2}=l_{2}\left( N\right) 
$ and $l_{2}^{s}=l_{2}^{s}\left( N\right) $ (see $\left[ \text{23, \S\ 1.18}%
\right] $). Let $A$ be the operator in $l_{2}\left( N\right) $ defined by%
\[
\text{ }D\left( A\right) =\left\{ u=\left\{ u_{j}\right\} \in l_{2},\text{ }%
\left\Vert u\right\Vert _{l_{2}^{s}\left( N\right) }=\left(
\sum\limits_{j=1}^{N}2^{sj}\left\vert u_{j}\right\vert ^{2}\right) ^{\frac{1%
}{2}}<\infty \right\} , 
\]

\[
A=\left[ a_{mj}\right] \text{, }a_{mj}=G_{m}\left( x\right) 2^{sj},\text{ }%
s>0,\text{ }m,j=1,2,...,N,\text{ }N\in \mathbb{N} 
\]

\bigskip and 
\[
\text{ }D\left( V\left( x,t\right) \right) =l_{2}, 
\]

\[
V\left( x,t\right) =\left[ b_{mj}\left( x,t\right) \right] \text{, }%
b_{mj}\left( x,t\right) =g_{m}\left( x,t\right) 2^{sj},\text{ }%
m,j=1,2,...,N. 
\]

Let 
\[
X_{2}=L^{2}\left( R^{n};l_{2}\right) ,\text{ }X_{2}\left( A\right)
=L^{2}\left( R^{n};l_{2}^{s}\right) ,\text{ }Y^{k,2}=W^{k,2}\left(
R^{n};l_{2}\right) . 
\]

\ From Theorem 1 we obtain the following result

\textbf{Theorem 6.1. }Suppose $a_{mj}$ are bounded cotinious on $R^{n}$ and $%
b_{mj}$\ \ are bounded functions on $R^{n}\times \left[ 0,T\right] .$ Assume
there exist the constants $a_{0},$ $a_{1},$ $a_{2}>0$ such that for any $%
k\in \mathbb{Z}^{+}$ a solution $u\in C\left( \left[ 0,1\right] ;X_{2}\left(
A\right) \right) $ of $\left( 6.1\right) $ satisfy 
\[
\dint\limits_{R^{n}}\left\Vert u\left( x,0\right) \right\Vert
_{l_{2}}^{2}e^{2a_{0}\left\vert x\right\vert ^{p}}dx<\infty ,\text{ for }%
p\in \left( 1,2\right) , 
\]%
\[
\dint\limits_{R^{n}}\left\Vert u\left( x,1\right) \right\Vert
_{l_{2}}^{2}e^{2k\left\vert x\right\vert ^{p}}dx<a_{2}e^{2a_{1}k^{\frac{q}{%
q-p}}},\text{ }\frac{1}{p}+\frac{1}{q}=1. 
\]%
Moreover, there exists $M_{p}>0$ such that

\[
a_{0}a_{1}^{p-2}>M_{p}. 
\]

Then $u\left( x,t\right) \equiv 0.$

\ \textbf{Proof.} It is easy to see that $A$ is a symmetric operator in $%
l_{2}$ and other conditions of Theorem 1 are satisfied. Hence, from Teorem 1
we obtain the conculision.

\begin{center}
\textbf{7. Unique continuation properties for nonlinear anisotropic
Schredinger equation }

\ \ \ \ \ \ \ \ \ \ \ \ \ \ \ \ \ \ \ \ \ \ \ \ \ \ \ \ \ \ \ \ \ \ \ \ \ \
\ \ \ \ \ \ \ \ \ \ 
\end{center}

The regularity property of BVPs for elliptic equations\ were studied e.g. in 
$\left[ \text{1, 2}\right] $. Let $\Omega =R^{n}\times G$, $G\subset R^{d},$ 
$d\geq 2$ is a bounded domain with $\left( d-1\right) $-dimensional boundary 
$\partial G$. Let us consider the following problem

\begin{equation}
i\partial _{t}u+\Delta _{x}u+\sum\limits_{\left\vert \alpha \right\vert \leq
2m}a_{\alpha }\left( x,y\right) D_{y}^{\alpha }u\left( x,y,t\right) +F\left(
u,\bar{u}\right) u=0,\text{ }  \tag{7.1}
\end{equation}%
\[
\text{ }x\in R^{n},\text{ }y\in G,\text{ }t\in \left[ 0,1\right] ,\text{ } 
\]

\begin{equation}
B_{j}u=\sum\limits_{\left\vert \beta \right\vert \leq m_{j}}\ b_{j\beta
}\left( y\right) D_{y}^{\beta }u\left( x,y,t\right) =0\text{, }x\in R^{n},%
\text{ }y\in \partial G,\text{ }j=1,2,...,m,  \tag{7.2}
\end{equation}%
where $a_{\alpha },$ $b_{j\beta }$ are the complex valued functions, $\alpha
=\left( \alpha _{1},\alpha _{2},...,\alpha _{d}\right) $, $\beta =\left(
\beta _{1},\beta _{2},...,\beta _{d}\right) ,$ $\mu _{i}<2m,$ $K=K\left(
x,y,t\right) $ and 
\[
D_{x}^{k}=\frac{\partial ^{k}}{\partial x^{k}},\text{ }D_{j}=-i\frac{%
\partial }{\partial y_{j}},\text{ }D_{y}=\left( D_{1,}...,D_{d}\right) ,%
\text{ }y=\left( y_{1},...,y_{d}\right) . 
\]

$\ $

\bigskip Let%
\[
\xi ^{\prime }=\left( \xi _{1},\xi _{2},...,\xi _{d-1}\right) \in R^{d-1},%
\text{ }\alpha ^{\prime }=\left( \alpha _{1},\alpha _{2},...,\alpha
_{d-1}\right) \in Z^{d-1},\text{ } 
\]%
\[
\text{ }A\left( y_{0},\xi ^{\prime },D_{y}\right) =\sum\limits_{\left\vert
\alpha ^{\prime }\right\vert +j\leq 2m}a_{\alpha ^{\prime }}\left(
y_{0}\right) \xi _{1}^{\alpha _{1}}\xi _{2}^{\alpha _{2}}...\xi
_{d-1}^{\alpha _{d-1}}D_{y}^{j}\text{ for }y_{0}\in \bar{G} 
\]%
\[
B_{j}\left( y_{0},\xi ^{\prime },D_{y}\right) =\sum\limits_{\left\vert \beta
^{\prime }\right\vert +j\leq m_{j}}b_{j\beta ^{\prime }}\left( y_{0}\right)
\xi _{1}^{\beta _{1}}\xi _{2}^{\beta _{2}}...\xi _{d-1}^{\beta
_{d-1}}D_{y}^{j}\text{ for }y_{0}\in \partial G. 
\]

Let 
\[
X_{2}=L^{2}\left( R^{n};L^{2}\left( G\right) \right) =L^{2}\left( \Omega
\right) ,\text{ }X_{2}\left( A\right) =L^{2}\left( R^{n};W^{2m,2}\left(
G\right) \right) ,\text{ }Y^{k,2}=W^{k,2}\left( R^{n};L^{2}\left( G\right)
\right) . 
\]

\textbf{Theorem 7.1}. Let the following conditions be satisfied:

\bigskip (1) $G\in C^{2}$, $a_{\alpha }\in C\left( \bar{\Omega}\right) $ for
each $\left\vert \alpha \right\vert =2m$ and $a_{\alpha }\in L_{\infty
}\left( G\right) $ for each $\left\vert \alpha \right\vert <2m$;

(2) $b_{j\beta }\in C^{2m-m_{j}}\left( \partial G\right) $ for each $j$, $%
\beta $ and $\ m_{j}<2m$, $\sum\limits_{j=1}^{m}b_{j\beta }\left( y^{\prime
}\right) \sigma _{j}\neq 0,$ for $\left\vert \beta \right\vert =m_{j},$ $%
y^{^{\shortmid }}\in \partial G,$ where $\sigma =\left( \sigma _{1},\sigma
_{2},...,\sigma _{d}\right) \in R^{d}$ is a normal to $\partial G$ $;$

(3) for $y\in \bar{G}$, $\xi \in R^{d}$, $\lambda $ with $\left\vert \arg
\lambda \right\vert \leq \varphi $ for $0\leq \varphi <\pi $, $\left\vert
\xi \right\vert +\left\vert \lambda \right\vert \neq 0$ let $\lambda +$ $%
\sum\limits_{\left\vert \alpha \right\vert =2m}a_{\alpha }\left( y\right)
\xi ^{\alpha }\neq 0$;

(4) for each $y_{0}\in \partial G$ local BVP in local coordinates
corresponding to $y_{0}$:%
\[
\lambda +A\left( y_{0},\xi ^{\prime },D_{y}\right) \vartheta \left( y\right)
=0, 
\]

\[
B_{j}\left( y_{0},\xi ^{\prime },D_{y}\right) \vartheta \left( 0\right)
=h_{j}\text{, }j=1,2,...,m 
\]%
has a unique solution $\vartheta \in C_{0}\left( \mathbb{R}_{+}\right) $ for
all $h=\left( h_{1},h_{2},...,h_{d}\right) \in \mathbb{C}^{d}$ and for $\xi
^{\prime }\in R^{d-1};$

(5) there exist positive constants $b_{0}$ and $\theta $ such that a
solution $u\in C\left( \left[ -1,1\right] ;X_{2}\left( A\right) \right) $ of 
$\left( 7.1\right) -\left( 7.2\right) $ satisfied 
\[
\left\Vert F\left( u,\bar{u}\right) \right\Vert _{L^{2}\left( G\right) }\leq
b_{0}\left\Vert u\right\Vert _{L^{2}\left( G\right) }^{\theta }\text{ for }%
\left\Vert u\right\Vert _{X_{2}\left( A\right) }>1; 
\]

(6) Suppose 
\[
\left\Vert u\left( .,t\right) \right\Vert _{L^{2}\left( \Omega \right)
}=\left\Vert u\left( .,0\right) \right\Vert _{L^{2}\left( \Omega \right)
}=\left\Vert u_{0}\right\Vert _{L^{2}\left( \Omega \right) }=a 
\]%
for $t\in \left[ -1,1\right] $ and that $\left( 2.15\right) $ holds with $%
Q\left( .\right) $ satisfies $\left( 2.16\right) $ for $H=L^{2}\left(
G\right) .$

\ If $p>p\left( \theta \right) =\frac{2\left( \theta n-2\right) }{\left(
\theta n-1\right) },$ then $a\equiv 0.$

\ \textbf{Proof. }Let us consider operators $A$ and $V\left( x,t\right) $ in 
$H=L^{2}\left( G\right) $ that are defined by the equalities 
\[
D\left( A\right) =\left\{ u\in W^{2m,2}\left( G\right) \text{, }B_{j}u=0,%
\text{ }j=1,2,...,m\text{ }\right\} ,\ Au=\sum\limits_{\left\vert \alpha
\right\vert \leq 2m}a_{\alpha }\left( y\right) D_{y}^{\alpha }u\left(
y\right) . 
\]

Then the problem $\left( 7.1\right) -\left( 7.2\right) $ can be rewritten as
the problem $\left( 2.12\right) $, where $u\left( x\right) =u\left(
x,.\right) ,$ $f\left( x\right) =f\left( x,.\right) $,\ $x\in R^{n}$ are the
functions with values in\ $H=L^{2}\left( G\right) $. By virtue of $\left[ 
\text{1}\right] $ operator $A+\mu $ is positive in $L^{2}\left( G\right) $
for sufficiently large $\mu >0$. Moreover, in view of (1)-(6) all conditons
of Theorem 3 are hold. Then Theorem 3 implies the assertion.

\begin{center}
\textbf{8.} \textbf{The Wentzell-Robin type mixed problem for Schredinger
equations}
\end{center}

Consider the problem $\left( 1.5\right) -\left( 1.6\right) $. \ Let 
\[
\sigma =R^{n}\times \left( 0,1\right) \text{, }X_{2}=L^{2}\left( \sigma
\right) ,\text{ }Y^{2,k}=W^{2,k}\left( \sigma \right) . 
\]

Suppose $\nu =\left( \nu _{1},\nu _{2},...,\nu _{n}\right) $ are nonnegative
real numbers. In this section, from Theorem 1 we obtain the following result:

\bigskip \textbf{Theorem 8.1. } Suppose the following conditions are
satisfied:

(1)\ $a$ \ and $b$ are a real-valued functions on $\left( 0,1\right) $ and $%
a\left( t\right) >0$. Moreover$,$ $a\left( .\right) $ is bounded continious
function on $R^{n}\times \left[ 0,1\right] $ and%
\[
\exp \left( -\dint\limits_{\frac{1}{2}}^{x}b\left( t\right) a^{-1}\left(
t\right) dt\right) \in L_{1}\left( 0,1\right) ; 
\]

(2) $u_{1},$ $u_{2}\in C\left( \left[ 0,1\right] ;Y^{2,k}\right) $ are
strong solutions of $(1.2)$ with $k>\frac{n}{2};$

(3) $F:\mathbb{C}\times \mathbb{C}\rightarrow \mathbb{C},$ $F\in C^{k}$, $%
F\left( 0\right) =\partial _{u}F\left( 0\right) =\partial _{\bar{u}}F\left(
0\right) =0;$

(4) there exist positive constants $\alpha $ and $\beta $ such that 
\[
e^{\frac{\left\vert \alpha x\right\vert ^{p}}{p}}\left( u_{1}\left(
.,0\right) -u_{2}\left( .,0\right) \right) \in X_{2},\text{ }e^{^{\frac{%
\left\vert \beta x\right\vert ^{q}}{q}}}\left( u_{1}\left( .,0\right)
-u_{2}\left( .,0\right) \right) \in X_{2}, 
\]

with 
\[
p\in \left( 1,2\right) \text{, }\frac{1}{p}+\frac{1}{q}=1; 
\]

(5) there exists $N_{p}>0$ such that 
\[
\alpha \beta >N_{p}.\text{ } 
\]

Then $u_{1}\left( x,t\right) \equiv u_{2}\left( x,t\right) .$

\ \textbf{Proof.} Let $H=L^{2}\left( 0,1\right) $ and $A$ is a operator
defined by $\left( 1.4\right) .$ Then the problem $\left( 1.5\right) -\left(
1.6\right) $ can be rewritten as the problem $\left( 1.2\right) $. By virtue
of $\left[ \text{10, 11}\right] $ the operator $A$ generates analytic
semigroup in $L^{2}\left( 0,1\right) $. Hence, by virtue of (1)-(5) all
conditons of Theorem 2 are satisfied. Then Theorem 2 implies the assertion.

\begin{center}
\textbf{Acknowledgements}
\end{center}

The author would like to express a gratitude to Dr. Neil. Course for his
useful advice in English in preparing of this paper

\textbf{References}\ \ 

\begin{enumerate}
\item H. Amann, Linear and quasi-linear equations,1, Birkhauser, Basel
(1995).

\item A. Benedek, A. Calder\`{o}n, R. Panzone, Convolution operators on
Banach space valued functions, Proc. Nat. Acad. Sci. USA, 48(1962), 356--365

\item A. Bonami, B. Demange, A survey on uncertainty principles related to
quadratic forms, Collect. Math., V. Extra, (2006), 1--36.

\item A.-P. Calder\'{o}n, Commutators of singular integral operators, Proc.
Nat. Acad. Sci. U.S.A. 53(1965), 1092--1099.

\item C. E. Kenig, G. Ponce, L. Vega, On unique continuation for nonlinear
Schr%
%TCIMACRO{\U{a8}}%
%BeginExpansion
\"{}%
%EndExpansion
odinger equations, Comm. Pure Appl. Math. 60 (2002), 1247--1262.

\item L. Escauriaza, C. E. Kenig, G. Ponce, L. Vega, On uniqueness
properties of solutions of Schr\"{o}dinger Equations, Comm. PDE. 31, 12
(2006) 1811--1823.

\item L. Escauriaza, C. E. Kenig, G. Ponce, and L. Vega, Uncertainty
principle of Morgan type and Schr\"{o}dinger evolution, J. London Math. Soc.
81, (2011) 187--207.

\item L. Escauriaza, C. E. Kenig, G. Ponce, and L. Vega, Hardy's uncertainty
principle, convexity and Schr\"{o}dinger Evolutions, J. European Math. Soc.
10, 4 (2008) 883--907.

\item L. H\"{o}rmander, A uniqueness theorem of Beurling for Fourier
transform pairs, Ark. Mat. 29, 2 (1991) 237--240.

\item J. A. Goldstain, Semigroups of Linear Operators and Applications,
Oxford University Press, Oxfard (1985).

\item A. Favini, G. R. Goldstein, J. A. Goldstein and S. Romanelli,
Degenerate Second Order Differential Operators Generating Analytic
Semigroups in $L_{p}$ and $W^{1,p}$, Math. Nachr. 238 (2002), 78 --102.

\item V. Keyantuo, M. Warma, The wave equation with Wentzell--Robin boundary
conditions on Lp-spaces, J. Differential Equations 229 (2006) 680--697.

\item S. G. Krein, Linear Differential Equations in Banach space, American
Mathematical Society, Providence, (1971).

\item A. Lunardi, Analytic Semigroups and Optimal Regularity in Parabolic
Problems, Birkhauser (2003).

\item J-L. Lions, E. Magenes, Nonhomogenous Boundary Value Broblems, Mir,
Moscow (1971).

\item V. B. Shakhmurov, Imbedding theorems and their applications to
degenerate equations, Differential equations, 24 (4) (1988), 475-482.

\item V. B. Shakhmurov, Linear and nonlinear abstract equations with
parameters, Nonlinear Anal-Theor., 73 (2010), 2383-2397.

\item R. Shahmurov, On strong solutions of a Robin problem modeling heat
conduction in materials with corroded boundary, Nonlinear Anal., Real World
Appl., 13, (1)( 2011), 441-451.

\item R. Shahmurov, Solution of the Dirichlet and Neumann problems for a
modified Helmholtz equation in Besov spaces on an annuals, J. Differential
equations, 249(3)(2010), 526-550.

\item C. Segovia, J. L.Torrea, Vector-valued commutators and applications,
Indiana Univ. Math.J. 38(4) (1989), 959--971.

\item E. M. Stein, R. Shakarchi, Princeton, Lecture in Analysis II. Complex
Analysis, Princeton University Press (2003).

\item A. Stefanov, Strichartz estimates for the Schrodinger equation with
radial data, Proc. Amer. Math. Soc. 129 (2001), 1395-1401.

\item B. Simon, M. Schechter, Unique Continuation for Schrodinger Operators
with unbounded potentials, J. Math. Anal. Appl., 77 (1980), 482-492.

\item H. Triebel, Interpolation theory, Function spaces, Differential
operators, North-Holland, Amsterdam (1978).

\item S. Yakubov and Ya. Yakubov, Differential-operator Equations. Ordinary
and Partial \ Differential Equations, Chapman and Hall /CRC, Boca Raton
(2000).

\item Y. Xiao and Z. Xin, On the vanishing viscosity limit for the 3D
Navier-Stokes equations with a slip boundary condition. Comm. Pure Appl.
Math. 60 (7) (2007), 1027--1055.
\end{enumerate}

\end{document}